\documentclass{article}
\usepackage{fullpage,amssymb,amsmath,amsthm}
\usepackage{epsfig}

\usepackage{color}
\definecolor{light-blue}{rgb}{0.8,0.85,1}
\definecolor{blue}{rgb}{0,0,1}
\definecolor{red}{rgb}{1,0,0}
\definecolor{green}{rgb}{0,1,0}

\begin{document}
\providecommand{\del}{{\Delta}} \providecommand{\nab}{{\nabla}}
\newtheorem{rhp}{Riemann-Hilbert Problem}
\newtheorem{prop}{Proposition} \newcommand{\cx}{{\mathbb C}}
\newtheorem{theorem}{Theorem} \newtheorem{corollary}{Corollary}
\newtheorem{Lemma}{Lemma} \newtheorem{definition}{Definition}
\newenvironment{remark}{$\triangleleft$ {\bf
Remark:}}{$\triangleright$} \newcommand{\nat}{{\mathbb N}}
\newcommand{\Z}{{\mathbb Z}} \newcommand{\mat}[1]{{\bf #1}}
\newcommand{\C}{{\mathbb C}}
\newcommand{\R}{{\mathbb R}}

\author{J. Baik\thanks{
Department of Mathematics, Princeton University and Department of
Mathematics, University of Michigan at Ann Arbor. Email:
\textbf{jbaik@math.princeton.edu}.}\and T.
Kriecherbauer\thanks{Fakult\"at f\"ur Mathematik,
Ruhr-Universit\"at Bochum. Email:
\textbf{Thomas.Kriecherbauer@ruhr-uni-bochum.de}.}
 \and K. T.-R. McLaughlin\thanks{Department of Mathematics, University
of North Carolina at Chapel Hill. Email:
\textbf{mcl@amath.unc.edu}.} \and P. D. Miller\thanks{Department
of Mathematics, University of Michigan. Email:
\textbf{millerpd@umich.edu}.}}

\title{Uniform Asymptotics for Polynomials Orthogonal
With Respect to a General Class of Discrete Weights and
Universality Results for Associated Ensembles: Announcement of
Results}

\date{\today} \maketitle

\begin{abstract}
We compute the pointwise asymptotics of orthogonal polynomials
with respect to a general class of pure point measures supported
on finite sets as both the number of nodes of the measure and also
the degree of the orthogonal polynomials become large. The class
of orthogonal polynomials we consider includes as special cases
the Krawtchouk and Hahn classical discrete orthogonal polynomials,
but is far more general. In particular, we consider nodes that are
not necessarily equally spaced.  The asymptotic results are given
with error bound for all points in the complex plane except for a
finite union of discs of arbitrarily small but fixed radii.  These
exceptional discs are the neighborhoods of the so-called band
edges of the associated equilibrium measure.  As applications, we
prove universality results for correlation functions of a general
class of discrete orthogonal polynomial ensembles, and in
particular we deduce asymptotic formulae with error bound for
certain statistics relevant in the random tiling of a hexagon with
rhombus-shaped tiles.

The discrete orthogonal polynomials are characterized in terms of
a a Riemann-Hilbert problem formulated for a meromorphic matrix
with certain pole conditions.  By extending the methods of [17,
22], we suggest a general and unifying approach to handle
Riemann-Hilbert problems in the situation when poles of the
unknown matrix are accumulating on some set in the asymptotic
limit of interest.
\end{abstract}


\section{Introduction}
This announcement concerns asymptotic properties of polynomials
that are orthogonal with respect to pure point measures supported
on finite sets.  Let $N\in\nat$ be fixed, and consider $N$
distinct real nodes $x_{N,0}<x_{N,1}<\dots <x_{N,N-1}$ to be
given; together the nodes make up the support of the pure point
measures we consider. We use the notation $X_N:=
\{x_{N,n}\}_{n=0}^{N-1}$ for the support set.  Along with nodes we
are given positive weights $w_{N,0},w_{N,1},\dots, w_{N,N-1}$,
which are the magnitudes of the point masses located at the
corresponding nodes. The {\em discrete orthogonal polynomials}
associated with this data are polynomials
$\{p_{N,k}(z)\}_{k=0}^{N-1}$ where $p_{N,k}(z)$ is of degree
exactly $k$ with a positive leading coefficient and where
\begin{equation}
\sum_{j=0}^{N-1}p_{N,k}(x_{N,j})p_{N,l}(x_{N,j})w_{N,j}=\delta_{kl}\,.
\label{eq:ortho}
\end{equation}
If $p_{N,k}(z)=c_{N,k}^{(k)}z^k + \dots + c_{N,k}^{(0)}$, then we
denote by $\pi_{N,k}(z)$ the associated monic polynomial
$p_{N,k}(z)/c_{N,k}^{(k)}$.  These polynomials exist and are
uniquely determined by the orthogonality conditions because the
inner product associated with (\ref{eq:ortho}) is positive
definite on ${\rm span}(1,z,z^2,\dots,z^{N-1})$ but is degenerate
on larger spaces of polynomials.  The polynomials $p_{N,k}(z)$ may
be constructed from the monomials by a Gram-Schmidt process.  A
general reference for properties of orthogonal polynomials
specific to the discrete case is the book of Nikiforov, Suslov,
and Uvarov \cite{NikiforovSU91}. In contrast to the discrete
orthogonal polynomials, we refer to the orthogonal polynomials
with respect to an absolutely continuous measure as the
\emph{continuous orthogonal polynomials}.

We use the notation $\mathbb{Z}_N$ for the set $\{0,1,2,\dots,N-1\}$.  Examples of classical discrete weights are \cite{AbramowitzS65}
\begin{itemize}
\item Krawtchouk weight: on the nodes $x_{N,j}:=(j+1/2)/N$ for
$j\in\mathbb{Z}_N$
in the interval $(0,1)$, we define the weight
\begin{equation}
  w^{\rm Kraw}_{N,j}(p,q):=
\frac{N^{N-1}\sqrt{pq}}{q^N\Gamma(N)} \binom{N-1}{j} p^jq^{N-1-j}\,.
\end{equation}
\item Hahn weight: on the infinite set of nodes $x_{N,j}:=(j+1/2)/N$
for
$j=0,1,2,\dots$, we define the weight
\begin{equation}
w_{N,j}(b,c,d):=
\frac{N^{N-1}}{\Gamma(N)}\cdot\frac{\Gamma(b)\Gamma(c+j)\Gamma(d+j)}
{\Gamma(j+1)\Gamma(b+j)\Gamma(c)\Gamma(d)}\,, \label{eq:Hahnraw}
\end{equation}
where $b$, $c$, and $d$ are real parameters. Special cases are the following.
First, taking $d=\alpha$ and $b=2-N-\beta$, with $\alpha, \beta >0$, and
taking the limit $c\to 1-N$, we obtain the weight
\begin{equation}\label{eq;Hahn1}
w_{N,j}^{\rm Hahn}(\alpha,\beta):=\frac{N^{N-1}}{\Gamma(N)}\cdot
\frac{\displaystyle\binom{j+\alpha-1}{j}\binom{N+\beta-2-j}{N-1-j}}
{\displaystyle\binom{N+\beta-2}{\beta-1}}\,,\hspace{0.2 in}
\mbox{for}\hspace{0.2 in}j\in{\mathbb Z}_N\,,
\end{equation}
which is defined on the \emph{finite} set of nodes $X_N:=(\mathbb{Z}_N+1/2)/N$. Second, taking $d=2-N-\beta$
and $b=\alpha$, with $\alpha, \beta>0$, and taking the limit $c\to 1-N$, we obtain
the
weight
\begin{equation}\label{eq;aHahn1}
w^{\rm Assoc}_{N,j}(\alpha,\beta):=\frac{N^{N-1}}{\Gamma(N)}
\cdot\frac{\Gamma(N)\Gamma(N+\beta-1)\Gamma(\alpha)}
{\Gamma(j+1)\Gamma(\alpha+j)\Gamma(N-j)\Gamma(N+\beta-1-j)}\,,\hspace{0.2
in} \mbox{for}\hspace{0.2 in}j\in{\mathbb Z}_N\,,
\end{equation}
which is again defined on the nodes $X_N:=(\mathbb{Z}_N+1/2)/N$.
\end{itemize}
\begin{remark}
 In \cite{AbramowitzS65}, the polynomials orthogonal
with respect to the weight \eqref{eq:Hahnraw} are called the Hahn polynomials, and
in \cite{BorodinO02}, the same weight is
called the Askey-Lesky weight.
For the two special cases where the weight is supported on a finite set of nodes, we now adopt the terminology used by Johansson \cite{Johansson00}, and thus refer to the weight (\ref{eq;Hahn1}) simply as the \emph{Hahn weight} (corresponding to  the Hahn polynomials) and we refer to the weight (\ref{eq;aHahn1}) as the
\emph{associated Hahn weight} (corresponding to the associated Hahn polynomials).
\end{remark}

Our goal is to establish the asymptotic behavior of the polynomials
$p_{N,k}(z)$ or their monic counterparts $\pi_{N,k}(z)$ in the limit
of large degree, assuming certain asymptotic properties of the nodes
and the weights.  In particular, the number of nodes must necessarily
increase to admit polynomials with arbitrarily large degree, and the
weights we consider involve an exponential factor with exponent
proportional to the number of nodes (such weights are sometimes called
{\em varying weights}).  We will obtain pointwise asymptotics with
precise error bound uniformly valid in the whole complex plane with
the exception of certain arbitrarily small open discs.  Our
assumptions on the nodes and weights include as special cases all
relevant classical discrete orthogonal polynomials, but are
significantly more general; in particular, we will consider nodes that
are not necessarily equally spaced.

\subsection{Motivation and applications.}

In the context of approximation theory, there has been recent
activity \cite{DKMVZstrong,DeiftKMVZ99} in the study of
polynomials orthogonal on the real axis with respect to general
continuous varying weights and the corresponding large degree
pointwise asymptotics. One of important application of the results
of \cite{DKMVZstrong, DeiftKMVZ99} is the proof of several
universality conjectures of random matrix theory. Thus, a natural
question to ask is whether it is possible to extend the results of
\cite{DKMVZstrong, DeiftKMVZ99} to handle discrete weights, and
obtain similar universality results for the so-called
\emph{discrete orthogonal polynomial ensembles} (see \S~\ref{sec;DOPensemble}).

Indeed, it has turned out recently that various problems of percolation
models, random tiling, queueing theory, non-intersecting paths and
representation theory can be reformulated as asymptotic questions
of discrete orthogonal polynomial ensembles with very concrete
weights (see for example, \cite{kurtj:shape, Johansson01,
BorodinO02} ). In these ensembles, the weights are all classical
(Meixner, Charlier, Krawtchouk or Hahn). Using integral
formulae for the corresponding orthogonal polynomials, the relevant
asymptotics have been analyzed except for the Hahn weight case.
The weights handled in this paper include the Hahn weight (and
also Krawtchouk), and hence as a corollary, we obtain new
asymptotic  results for the Hahn polynomials. As the discrete
orthogonal polynomial ensemble for the Hahn weight arises in the statistical analysis of  random
rhombus tilings of a hexagon, our asymptotic results for Hahn
polynomials yield new results on this problem (see
\S~\ref{sec;tiling}).

\subsection{Methodology.}\label{sec:method}

The method we use is the Riemann-Hilbert characterization of
discrete orthogonal polynomials, and an adaption of the Deift-Zhou
method for the steepest-descent analysis of Riemann-Hilbert
problems. The Riemann-Hilbert problem for discrete orthogonal
polynomials has poles instead of usual jump conditions on a
continuous contour, and the poles are accumulating in the limit of interest to form a continuum.

There has been some recent progress \cite{KamvissisMM02, Miller02}
in the integrable systems literature concerning the problem of
computing asymptotics for solutions of integrable nonlinear
partial differential equations ({\em e.g.} the nonlinear
Schr\"odinger equation) in the limit where the spectral data
associated with the solution via the inverse-scattering transform
is made up of a large number of discrete eigenvalues.
Significantly, inverse-scattering theory also exploits much of the
theory of matrix Riemann-Hilbert problems, and it turns out that
the discrete eigenvalues appear as poles in the corresponding
matrix-valued unknown.  So, the methods recently developed in the
context of inverse-scattering actually suggest a general scheme by
means of which an accumulation of poles in the matrix unknown can
be analyzed.

In this paper, we extend the method of \cite{KamvissisMM02,
Miller02} and suggest a general and unifying approach to handle
Riemann-Hilbert problems for the situation when poles are
accumulating. Especially we overcome the following two issues:
\begin{itemize}
\item[(a)] How to transform a Riemann-Hilbert problem with pole
conditions to a Riemann-Hilbert problem with an analytic jump condition on a
continuous contour so that a formal continuum limit of poles can
be rigorously justified and the Deift-Zhou method can be
applied.
\item[(b)]  How to handle the upper constraint of the so-called
equilibrium measure, and thus to correctly formulate an appropriate $g$-function.
\end{itemize}
See \S~\ref{sec;proof} for more information about these ideas.
Full details will be given in the paper corresponding to this announcement.

\subsection{Basic assumptions.}

We state here precise assumptions on the nodes $X_N$ and weights
$\{w_{N,j}\}$.

\subsubsection{Conditions on the nodes.}
\label{sec:C1}
\begin{enumerate}
\item
The nodes lie in a bounded open interval $(a,b)$ and are
distributed with a density $\rho^0(x)$.
\item
The density function $\rho^0(x)$ is real analytic in a complex
neighborhood of the closed interval $[a,b]$, and satisfies:
\begin{equation}
\int_a^b\rho^0(x)\,dx=1\,, \label{eq:rho0prob}
\end{equation}
and
\begin{equation}
\rho^0(x)>0\hspace{0.2 in}\mbox{strictly, for all $x\in[a,b]$.}
\label{eq:rho0nonzero}
\end{equation}
\item
The nodes are defined precisely in terms of the density function
$\rho^0(x)$ by the quantization rule
\begin{equation}
\int_a^{x_{N,j}}\rho^0(x)\,dx = \frac{2j+1}{2N} \label{eq:BS}
\end{equation}
for $N\in\nat$ and $j\in{\mathbb Z}_N$.
\end{enumerate}
\subsubsection{Conditions on the weights.}
\label{sec:C2}
\begin{enumerate}
\item
Without loss of generality, we write the weights in the form
\begin{equation}
w_{N,j}=(-1)^{N-1-j}e^{-NV_N(x_{N,j})}\mathop{\prod_{n=0}^{N-1}}_{n\neq
j}
(x_{N,j}-x_{N,n})^{-1}=e^{-NV_N(x_{N,j})}\mathop{\prod_{n=0}^{N-1}}_{n\neq
j} |x_{N,j}-x_{N,n}|^{-1}\,, \label{eq:weightform}
\end{equation}
where the family of functions $\{V_N(x)\}$ is {\em apriori}
specified only at the nodes.
\item
We assume that for each sufficiently large $N$, $V_N(x)$ may be
taken to be a real analytic function defined in a neighborhood $G$
of the closed interval $[a,b]$, and that
\begin{equation}
V_N(x)=V(x)+\frac{\gamma}{N} + \frac{\eta_N(x)}{N^2}
\label{eq:VNexpand}
\end{equation}
where $V(x)$ is a fixed real analytic function defined in $G$,
$\gamma$ is a constant, and
\begin{equation}
\limsup_{N\rightarrow\infty}\sup_{z\in G}|\eta_N(z)|<\infty\,.
\label{eq:etaNcontrol}
\end{equation}
\end{enumerate}

\begin{remark}
In some applications it is desirable to generalize further by
allowing $\gamma$ in \eqref{eq:VNexpand} to be a real analytic
function in $G$ with $\gamma'(z)$ not identically zero.  It is
possible to take into account such variation, but for simplicity
we take $\gamma$ to be constant in this paper. For classical cases
of Hahn and Krawtchouk weights, $\gamma$ is indeed constant in scalings of interest.
\end{remark}

The familiar examples of classical discrete orthogonal polynomials
correspond to nodes that are equally spaced, say on $(a,b)=(0,1)$
(in which case we have $\rho^0(x)\equiv 1$).  In this special
case, the product factor on the right-hand side of
(\ref{eq:weightform}) becomes simply
\begin{equation}
\mathop{\prod_{n=0}^{N-1}}_{n\neq j}|x_{N,j}-x_{N,n}|^{-1}=
\frac{N^{N-1}}{j!(N-j-1)!}\,.
\end{equation}
By Stirling's formula, taking the continuum limit of this factor
(that is, considering $N\to\infty$ with $j/N\to x$) shows that in
these cases the formula \eqref{eq:weightform} leads to a
continuous weight on $(0,1)$ of the form
\begin{equation}
  w(x) = \biggl( \frac{e^{-V(x)}}{x^x(1-x)^{1-x}} \biggr)^N
\end{equation}
up to an overall multiplicative constant.

Our choice of the form (\ref{eq:weightform}) for the weights is
motivated by several specific examples of classical discrete
orthogonal polynomials.  The form (\ref{eq:weightform}) is
sufficiently general for us to carry out useful calculations
related to proofs of universality conjectures arising in statistical problems like the random
rhombus tiling of a hexagon.

\subsubsection{Conditions on the equilibrium measure.}

There is an additional assumption, which is not as explicit as the
previous two conditions. This assumption will be explained in the
next subsection.

\subsection{The equilibrium energy problem and third assumption on the weights.}
\label{sec:equilibrium}

It has been recognized for some time (see \cite{KuijlaarsR98} and
references therein) that, as in the continuous orthogonal
polynomial cases, the asymptotic behavior of discrete orthogonal
polynomials, in particular the distribution of zeros
in $(a,b)$, is related to a constrained equilibrium problem for
logarithmic potentials in a field $\varphi(x)$ given by the
formula
\begin{equation}
\varphi(x):=V(x)+\int_a^b\log|x-y|\rho^0(y)\,dy
\label{eq:fielddef}
\end{equation}
for $x\in (a,b)$.  We can also view $\varphi(x)$ as being defined
via a continuum limit:
\begin{equation}
\varphi(x)=-\lim_{N\rightarrow\infty}\frac{\log(w_{N,j})}{N}
\end{equation}
where $w_{N,j}$ is expressed in terms of $x_{N,j}$ which in turn
is identified with $x$. Thus at the moment we are working with the
\emph{formal continuum limit} of the weight
$w_{N,j}$.

In the specific context of this paper, the field $\varphi(x)$ is a
real analytic function in the open interval $(a,b)$ because $V(x)$
and $\rho^0(x)$ are real analytic functions in a neighborhood of
$[a,b]$. Unlike $V(x)$ and $\rho^0(x)$, however, the field
$\varphi(x)$ does not extend analytically beyond the endpoints of
$(a,b)$ due to the condition \eqref{eq:rho0nonzero}.

Given the parameter $c\in (0,1)$, which has the interpretation of the ratio of the degree $k$ of the polynomial of interest to the number $N$ of nodes,  and the field $\varphi(x)$ as above, consider the
quadratic functional
\begin{equation}
E_c[\mu]:=c\int_a^b\int_a^b\log\frac{1}{|x-y|} \,
d\mu(x)\,d\mu(y)+ \int_a^b \varphi(x)\,d\mu(x) \label{eq:energy}
\end{equation}
of Borel measures $\mu$ on $[a,b]$. Let $\mu^c_{\rm min}$ be the
measure that minimizes $E_c[\mu]$ over the class of measures
satisfying the upper and lower constraints
\begin{equation}
0\le \int_{x\in{\mathcal
B}}d\mu(x)\le\frac{1}{c}\int_{x\in{\mathcal B}}\rho^0(x)\,dx
\label{eq:constraints}
\end{equation}
for all Borel sets ${\mathcal B}\subset [a,b]$, and the
normalization condition
\begin{equation}
\int_a^bd\mu(x)=1\,. \label{eq:Lagrange}
\end{equation}
The existence of a unique minimizer under the conditions
enumerated in \S~\ref{sec:C1} and \S~\ref{sec:C2} follows from the
Gauss-Frostman Theorem; see \cite{ST} for details.  We will often
refer to the minimizer as the {\em equilibrium measure}.  It has been shown
\cite{KuijlaarsR98} that the equilibrium measure is the weak limit of the normalized
counting measure of the zeros of $p_{N,k}(z)$ in the limit $N\rightarrow\infty$ with
$c=k/N$ fixed.

That a variational problem plays a central role in asymptotic
behavior is a familiar theme in the theory of orthogonal
polynomials.  The key new feature contributed by discreteness is
the appearance of the upper constraint on the equilibrium measure
({\em i.e.} the upper bound in \eqref{eq:constraints}).  The upper
constraint can be traced to the following well-known fact.

\begin{prop}
Each discrete orthogonal polynomial $p_{N,k}(z)$ has $k$ simple
real zeros.  All zeros lie in the range $x_{N,0} < z < x_{N,N-1}$
and no more than one zero lies in the closed interval
$[x_{N,n},x_{N,n+1}]$ between any two consecutive nodes.
\label{prop:confine}
\end{prop}

Thus, the presence of the upper constraint proportional to the local density of nodes is necessary
for the interpretation of the equilibrium measure as the weak limit of the normalized
counting measure of zeros.

The theory of the ``doubly constrained'' variational problem we
are considering is well-established.  In particular, the analytic
properties we assume of $V(x)$ and $\rho^0(x)$ turn out to be
unnecessary for the mere existence of the minimizer.  However, it
has been shown \cite{Kuijlaars00} that analyticity of $V(x)$ and
$\rho^0(x)$ in a neighborhood of $[a,b]$ guarantees that $\mu_{\rm
min}^c$ is continuously differentiable with respect to $x\in
(a,b)$.  Moreover, the derivative $d\mu_{\rm min}^c/dx$ is
piecewise analytic, with a finite number of points of
nonanalyticity that may not occur at any $x$ where both (strict)
inequalities $d\mu_{\rm min}^c/dx(x)>0$ and $d\mu_{\rm
min}^c/dx(x)<\rho^0(x)/c$ hold.  We want to exploit these facts,
which is why we have chosen to restrict attention to analytic
functions $V(x)$ and $\rho^0(x)$.

For a method of computing the equilibrium measure from the coefficients
in the three-term recurrence relation for a special
class of discrete weights, see \cite{KuijlaarsR98} and reference therein. See also
\cite{DKM} for continuous weights.

For simplicity of exposition we want to exclude certain nongeneric
phenomena that may occur even under conditions of analyticity of
$V(x)$ and $\rho^0(x)$.  Therefore we introduce the following
assumptions.

\subsubsection{Third assumption; conditions on the equilibrium measure.}
\label{sec:C3} Let $\underline{\mathcal F}\subset [a,b]$ denote
the closed set of $x$-values where $d\mu_{\rm min}^c/dx(x)=0$.
Let $\overline{\mathcal F}\subset [a,b]$ denote the closed set of
$x$-values where $d\mu_{\rm min}^c/dx(x)=\rho^0(x)/c$.
\begin{enumerate}
\item
Each connected component of $\underline{\mathcal F}$ and
$\overline{\mathcal F}$ has a nonempty interior.  Therefore
$\underline{\mathcal F}$ and $\overline{\mathcal F}$ are both
finite unions of closed intervals, where each closed interval that is part of the union contains more than one
point.  Note:  this does not exclude the possibility that either $\underline{\mathcal F}$ or $\overline{\mathcal F}$ might be empty.
\item
For each open subinterval $I$ of $(a,b)\setminus
\underline{\mathcal F}\cup \overline{\mathcal F}$ and each limit
point $z_0\in\underline{\mathcal F}$ of $I$, we have
\begin{equation}
\lim_{x\rightarrow z_0, x\in I}\frac{1}{\sqrt{|x-z_0|}}
\frac{d\mu_{\rm min}^c}{dx}(x) = K\hspace{0.2 in}\mbox{with
$0<K<\infty$} \label{eq:rootlower}
\end{equation}
and for each limit point $z_0\in\overline{\mathcal F}$ of $I$, we
have
\begin{equation}
\lim_{x\rightarrow z_0, x\in I}\frac{1}{\sqrt{|x-z_0|}}
\left[\frac{1}{c}\rho^0(x)-\frac{d\mu_{\rm min}^c}{dx}(x)\right] =
K\hspace{0.2 in}\mbox{with $0<K<\infty$}\,. \label{eq:rootupper}
\end{equation}
Therefore the derivative of the minimizing measure meets each
constraint exactly like a square root.
\item
A constraint is active at each endpoint: $\{a,b\}\subset
\underline{\mathcal F}\cup\overline{\mathcal F}$.
\end{enumerate}
It is difficult to translate these conditions on $\mu_{\rm min}^c$
into sufficient conditions on $c$, $V(x)$, and $\rho^0(x)$.
However, there is a sense in which they are satisfied generically.

\begin{remark}
Relaxing the condition that a constraint should be active at each
endpoint requires specific local analysis near these two points.
We expect that a constraint being active at each endpoint is a
generic phenomenon in the sense that the opposite situation occurs
only for isolated values of $c$.  We know this statement to be
true in all relevant classical cases.  For the Krawtchouk
polynomials only the values $c=p$ or $c=q=1-p$ correspond to an
equilibrium measure that is not constrained at both endpoints (see
\cite{DragnevS00}).  The situation is similar for the Hahn
polynomials.
\end{remark}

\subsubsection{Voids, bands and saturated regions.}

Under the conditions enumerated in \S~\ref{sec:C1},
\S~\ref{sec:C2}, and \S~\ref{sec:C3}, the minimizer $\mu^c_{\rm
min}$ partitions $(a,b)$ into three kinds of subintervals, a
finite number of each, and each having a nonempty interior.  There
is a real constant $\ell_c$, the Lagrange multiplier associated
with the condition (\ref{eq:Lagrange}), so that with the
variational derivative defined as
\begin{equation}
\frac{\delta E_c}{\delta\mu}(x):=-2c\int_a^b\log|x-y|\,d\mu(y) +
\varphi(x)\,, \label{eq:variationalderivative}
\end{equation}
we have when $\mu=\mu^c_{\rm min}$ the following types of
subintervals.
\begin{definition}[Voids]
A void $\Gamma$ is an open subinterval of $[a,b]$ of maximal
length in which $\mu^c_{\rm min}(x)\equiv 0$, and thus the
minimizer realizes the lower constraint.  For $x\in\Gamma$ we have
the strict inequality
\begin{equation}
\frac{\delta E_c}{\delta\mu}(x)>\ell_c\,.
\label{eq:voidinequality}
\end{equation}
\end{definition}
\begin{definition}[Bands]
A band $I$ is an open subinterval of $[a,b]$ of maximal length
where $\mu^c_{\rm min}(x)$ is a measure with a real analytic
density satisfying $0<d\mu^c_{\rm min}/dx<\rho^0(x)/c$, and thus
variations of the minimizer are free.  For $x\in I$ we have the
equilibrium condition
\begin{equation}
\frac{\delta E_c}{\delta\mu}(x)\equiv\ell_c\,.
\label{eq:equilibrium}
\end{equation}
\end{definition}
\begin{definition}[Saturated regions]
A saturated region $\Gamma$ is an open subinterval of $[a,b]$ of
maximal length in which $d\mu^c_{\rm min}/dx\equiv \rho^0(x)/c$,
and thus the minimizer realizes the upper constraint.  For
$x\in\Gamma$ we have the strict inequality
\begin{equation}
\frac{\delta E_c}{\delta\mu}(x)<\ell_c\,.
\label{eq:saturatedregioninequality}
\end{equation}
\end{definition}

\begin{figure}[h]
\begin{center}
\input{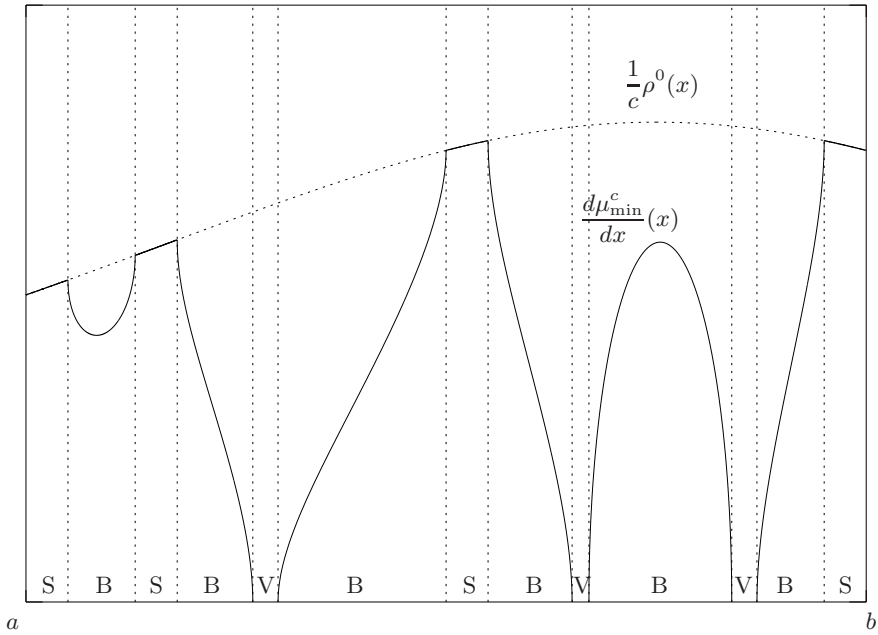}
\end{center}
\caption{\em The hypothetical minimizer illustrated here partitions the interval $(a,b)$ into voids (denoted V), bands (denoted B), and saturated regions (denoted S) as shown.}
\label{fig:intervals}
\end{figure}

See Figure~\ref{fig:intervals} for an example of the voids, bands and saturated regions associated with a hypothetical equilibrium measure.

Voids and saturated regions will also be called {\em gaps} when it
is not necessary to distinguish between these two types of
intervals. The closure of the union of all subintervals of the
three types defined above is the interval $[a,b]$.  From condition
1 in \S~\ref{sec:C3} above, bands cannot be adjacent to each
other; a band that is not adjacent to an endpoint of $[a,b]$ has
on each side either a void or a saturated region.

Some of our asymptotic results for the discrete orthogonal polynomials
under the above assumptions are stated in
\S~\ref{sec;results}. In \S~\ref{sec;applications}, the
asymptotics of discrete orthogonal polynomials are applied to
discrete orthogonal polynomial ensembles, and asymptotics of corresponding
correlation functions are thus obtained.  By specializing to the Hahn ensemble, we arrive at  specific results relevant in the problem of random rhombus tiling of a
hexagon. Our new methods for
the asymptotic analysis of general discrete orthogonal polynomials are
discussed in \S~\ref{sec;proof}.

\section{Results: Pointwise Asymptotics of Orthogonal
Polynomials}\label{sec;results}

As stated in the Introduction, we obtain pointwise asymptotics of
the orthogonal polynomials $p_{N,k}(z)$ for $z$ in the complex
plane except for a finite union of discs of arbitrarily small but
fixed radii in terms of a log transform of the equilibrium
measure. These discs are centered at the edges of the bands. The
difficulty at the edge of the bands will be explained below
(see \S~\ref{sec;proof}). We hope to be able to handle the band edge
problem in our future publication.

We actually present here formulae for the monic polynomials $\pi_{N,k}(z)$, rather
than the normalized polynomials $p_{N,k}(z)$.
The asymptotics are given by different formulae in five
different regions of the complex plane excluding small discs
around the edges of the bands. These regions are (recall that
$[a,b]$ is the interval where the nodes are accumulating)
\begin{itemize}
\item[(a)] Outside the interval $[a,b]$.
\item[(b)] Voids in $[a,b]$.
\item[(c)] Bands in $[a,b]$.
\item[(d)] Saturated regions in $[a,b]$.
\item[(e)] Near the endpoints of $[a,b]$ adjacent to a saturated
region.
\end{itemize}
The results for the first three cases are analogous to the
corresponding results for the continuous orthogonal polynomials
analyzed in \cite{DeiftKMVZ99}. Also the asymptotic formula for
$\pi_{N,k}$ near the endpoints of $[a,b]$ adjacent to a void is
analogous to the corresponding asymptotics of continuous
orthogonal polynomials whose weight is supported on a finite
interval.
The new
cases that did not occur in the continuous orthogonal polynomial
theory are the regions (d) and (e). In these regions, the discrete
nature of the support of the weights is strongly present. We here
present only the three regions (c), (d), (e). Asymptotic formulae
for $z$ in regions (a) and (b) will appear in the full version of this paper,
together with the asymptotics of the leading coefficient $c_{N,k}^{(k)}$ which completes the connection with the
polynomials $p_{N,k}(z)$.

For simplicity, we take the limit $k,N\to\infty$ while $c=k/N$ is
a fixed rational number. (Recall $k$ is the degree of the
orthogonal polynomial, and $N$ is the number of nodes.) This means
that if the fixed rational constant $c$ is represented in lowest
terms as $c=p/q$, then we are taking $k=Mp$ and $N=Mq$ for
$M\in\nat$. The more general case of $k/N = c +O(N^{-1})$ will be
considered in a future publication.

In the Theorems \ref{theorem:band} and \ref{theorem:upper} below,
the error bounds are different depending on the following two
situations of the support of the equilibrium measure.
\begin{itemize}
\item Case I: there is at least one void interval and at least
one saturated region.
\item Case II: there are only voids and bands, or only saturated
regions and bands.
\end{itemize}
We have obtained a better error bound for Case II than for
Case I. Whether this is only a technical point, or whether this is in the very
nature of discrete orthogonal polynomials is not clear yet.

\subsection{Bands.}

\begin{theorem}[Asymptotics of $\pi_{N,k}(z)$ in bands]
Assume the conditions enumerated in \S~\ref{sec:C1},
\S~\ref{sec:C2}, and \S~\ref{sec:C3}.  Let $c=k/N$ be fixed. Then,
uniformly for $z$ in any fixed compact subinterval in the interior
of a band $I$,
\begin{equation}
\pi_{N,k}(z)=\exp\left(k\int_a^b\log|z-x|\,d\mu_{\rm
min}^c(x)\right) \cdot\left[2A_I(z)\cos\left(k\pi\int_z^bd\mu_{\rm
min}^c(x) + \Phi_I(z)\right) +\varepsilon_N(z)\right]\,,
\end{equation}
where
\begin{equation}\label{eq;varepsilon}
\varepsilon_N(z)=\left\{\begin{array}{ll} \displaystyle
O\left(\frac{\log(N)}{N^{1/3}}\right)\,, &
\mbox{\rm Case I}\,,\\\\
\displaystyle O\left(\frac{\log(N)}{N^{2/3}}\right)\,, & \mbox{\rm
Case II}\,,
\end{array}\right.
\end{equation}
and $A_I(z)>0$ and $\Phi_I(z)$ are real functions defined in terms
of a Riemann theta function associated to the hyperelliptic
surface with cuts given by the bands. The functions $A_I$ and
$\Phi_I$ are uniformly bounded along with all derivatives.
\label{theorem:band}
\end{theorem}

\begin{remark}
The error estimates quoted above are derived from the asymptotic procedure we have
devised for the case of node distributions associated with a general analytic
density function $\rho^0(x)$.
In the special case of equally spaced nodes (when $\rho^0(x)$ is a constant function), we have recently found using a different procedure that the logarithmic term $\log(N)$
in the error bound $\varepsilon_N(z)$ can apparently be removed. At this time we do not know whether this alternate procedure can be modified for nonconstant $\rho^0(x)$ so as
to remove the logarithm from the estimates in all cases.
\end{remark}

\subsection{Saturated regions.}\label{sec:saturated}

\begin{theorem}[Asymptotics of $\pi_{N,k}(z)$ in saturated regions]
Assume the conditions enumerated in \S~\ref{sec:C1},
\S~\ref{sec:C2}, and \S~\ref{sec:C3}.  Let $c=k/N$ be fixed. Then,
uniformly for $z$ in a saturated region $\Gamma$, but bounded away
from any band edge points where the equilibrium measure becomes
unconstrained by a fixed distance and from the endpoints $a$ and
$b$ by a distance of size $N^{-2/3}$, we have
\begin{equation}
\begin{array}{rcl}
\displaystyle\pi_{N,k}(z)&=&\displaystyle
\exp\left(k\int_a^b\log|z-x|\,d\mu^c_{\rm min}(x)\right)\\\\
&&\displaystyle\hspace{0.2 in}\times\,\,
\left[\left(\phi_\Gamma(z)+\varepsilon_N(z) \right) \cdot 2
\cos\left(\pi N\int_z^b\rho^0(x)\,dx\right) + \mbox{\rm
exponentially small} \right]\,,
\end{array}
\label{eq:piNkupper}
\end{equation}
with $\varepsilon_N(z)$ having the same error bound as in
\eqref{eq;varepsilon},
where $\phi_\Gamma(z)$ is a real function defined in terms of a Riemann theta
function, uniformly bounded along with all
derivatives and having at most one zero in $\Gamma$. If
the saturated region $\Gamma$ is adjacent to an endpoint of
$[a,b]$, then $\phi_\Gamma(z)$ has no zeros in $\Gamma$.  The
exponentially small term is proportional to $\exp(N[\delta
E_c/d\mu-\ell_c])$ where the variational derivative is evaluated
on the equilibrium measure. \label{theorem:upper}
\end{theorem}

Since the zeros of the cosine function in (\ref{eq:piNkupper}) are
exactly the nodes of orthogonalization making up the set $X_N$,
and since the slope of the cosine is proportional to $N$ at the
nodes, we have the following.
\begin{corollary}[Exponential confinement of zeros]
Let $J$ be a closed subinterval of a saturated region where the
equilibrium measure achieves the upper constraint. Then the monic
discrete orthogonal polynomial $\pi_{N,k}(z)$ has a zero uniformly
close to each the nodes $x_{N,n}\in X_N\cap J$, with the possible
exception of one node.
 \label{corollary:exponential}
\end{corollary}

\begin{remark}
The factor $\phi_\Gamma(z)+\varepsilon_N(z)$ has at most one zero
in $\Gamma$, and the zero can be present for some $N$ and not
others. Also, the zero generally moves about in a quasiperiodic
manner as $N$ is varied. So it seems that one should regard the
situation in which this zero is exponentially close to one of the
nodes (which form a set of measure zero in $J$) as being anomalous
and quite rare. Therefore one should generally expect to see a
zero exponentially close to \emph{each} node in $X_N\cap J$.
\end{remark}

Let $K$ be a subinterval of the interval $J$ that is the subject
of Corollary~\ref{corollary:exponential} such that there is a
zero of $\pi_{N,k}(z)$ exponentially close to each node in
$X_N\cap K$ (so according to the previous remark one expects that
it is typically consistent to take $K=J$). The exponential
confinement of the zeros in $K$ has further consequences due to
the rigidity of the zeros for general discrete orthogonal
polynomials described in Proposition~\ref{prop:confine}.  A
particular zero $z_0\in K$ of $\pi_{N,k}(z)$, asymptotically
exponentially localized near a node $x_{N,n}$, can lie on one side
or the other of the node. But if $z_0$ lies to the right of
$x_{N,n}$, then it follows from Proposition~\ref{prop:confine}
that the smallest zero greater than $z_0$ must also lie to the
right of $x_{N,n+1}$ and so on, all the way to the right endpoint
of $K$.  Likewise, if $z_0$ lies to the left of $x_{N,n}$, then
all zeros in $K$ less than $z_0$ also lie to the left of the nodes
to which they are exponentially attracted.

\begin{figure}[h]
\begin{center}
\input{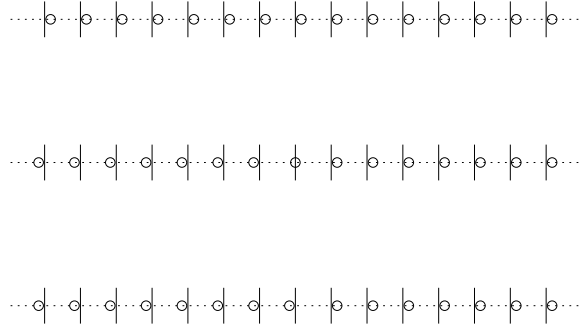}
\end{center}
\caption{\em The Hurwitz zeros are exponentially close to the
nodes of orthogonalization in saturated regions where the
equilibrium measure achieves its upper constraint.  Top: a pattern
without any dislocation, where all Hurwitz zeros (pictured as
circles) lie to the right of the nodes (vertical line segments) to
which they are exponentially attracted.  Bottom: a pattern with a
dislocation. Middle: there may only be one dislocation, but it can
move as parameters ({\em e.g.} $c$) are continuously varied and a
Hurwitz zero passes through one of the nodes.}
\label{fig:exponentialzeros}
\end{figure}

When we consider those zeros of $\pi_{N,k}(z)$ that converge
exponentially fast to the nodes (these are analogous to the {\em Hurwitz zeros} of the approximation theory
literature, whereas the possible lone zero of $\phi_\Gamma(z)$ would be called a {\em spurious zero}), we therefore see that there can be at most one
``dislocation'' ({\em i.e.} a closed interval of the form
$[x_{N,n},x_{N,n+1}]$ containing no Hurwitz zeros) in the pattern
of zeros lying to one side or the other of the nodes.  See
Figure~\ref{fig:exponentialzeros}.

\begin{remark}
It should perhaps be mentioned that there is nothing that prevents
a zero of $\pi_{N,k}(z)$ from coinciding {\em exactly} with one of
the nodes $x_{N,j}\in X_N$.
\end{remark}

Furthermore, due to Proposition~\ref{prop:confine}, a spurious
zero of $\pi_{N,k}(z)$ in the subinterval $K$ of the saturated
region $\Gamma$ can only occur if the pattern of Hurwitz zeros in
$K$ has a dislocation as in the bottom picture in Figure
\ref{fig:exponentialzeros}, in which case the spurious zero must
lie in the closed interval $[x_{N,n},x_{N,n+1}]$ associated with
the dislocation. Equivalently, the presence of a spurious zero in
$K\subset J\subset\Gamma$ indicates a dislocation in the pattern of
Hurwitz zeros.

>From the analysis we have presented it is not clear whether the
presence of a dislocation in the pattern of Hurwitz zeros implies
that the function $\phi_\Gamma(z)$ has a (spurious) zero in the
corresponding closed interval $[x_{N,n},x_{N,n+1}]$, or
equivalently whether the absence of any zeros of $\phi_\Gamma(z)$
in a saturated region $\Gamma$ means that the lone possible
dislocation in the pattern of Hurwitz zeros is indeed absent.

\subsection{Near hard edges adjacent to a saturated region.}

\begin{theorem}[Asymptotics of $\pi_{N,k}(z)$ near hard edges]
Assume the conditions enumerated in \S~\ref{sec:C1},
\S~\ref{sec:C2}, and \S~\ref{sec:C3}.  Let $c=k/N$ be fixed. If
the upper constraint is achieved at the endpoint $z=b$, then
uniformly for $b-CN^{-2/3}<z<b$,
\begin{equation}
\begin{array}{rcl}
\displaystyle \pi_{N,k}(z)&=&\displaystyle
\exp\left(k\int_a^b\log|z-x|\,d\mu^c_{\rm min}(x)\right)\\\\
&&\displaystyle\hspace{0.2 in}\times\,\,
\Bigg[\left(\phi_\Gamma(z)+O\left(\frac{\log(N)}{N^{1/3}}\right)
\right)\frac{\displaystyle\Gamma(1/2-\zeta_b)}
{\displaystyle\sqrt{2\pi}e^{\zeta_b}
(-\zeta_b)^{-\zeta_b}}2\cos\left(N\pi\int_z^b\rho^0(x)\,dx\right)
\\\\
&&\displaystyle\hspace{0.4 in}+\,\,\mbox{\rm exponentially
small}\Bigg]\,,
\end{array}
\end{equation}
as $N\rightarrow\infty$, and uniformly for $b<z<b+CN^{-2/3}$,
\begin{equation}
\pi_{N,k}(z)=\exp\left(k\int_a^b\log|z-x|\,d\mu_{\rm
min}^c(x)\right) \left[\phi^{\rm outside}(z)\frac{\displaystyle
\sqrt{2\pi}e^{-\zeta_b}\zeta_b^{\zeta_b}}
{\displaystyle\Gamma(1/2+\zeta_b)} +
O\left(\frac{\log(N)}{N^{1/3}} \right)\right]\,,
\end{equation}
as $N\rightarrow\infty$, where $\zeta_b:=N\rho^0(b)(z-b)$ and the
function $\phi_\Gamma(z)$ is same function that appears in
\eqref{eq:piNkupper}. The function $\phi^{\rm outside}(z)$ is real-valued, and nonvanishing, and like $\phi_\Gamma(z)$ is constructed from Riemann theta functions
and is along with all derivatives uniformly bounded in $z$ as $N\to\infty$.  There are similar formulae near the endpoint $z=a$ when the upper constraint is active there. The quantity
$\Gamma(1/2-\zeta_b)$ is the Euler gamma function, while the subscript
$\Gamma$ refers to the saturated region adjacent to
the hard edge. The exponentially small term is proportional to
$\exp(N[\delta E_c/\delta\mu-\ell_c])$ evaluated on the
equilibrium measure. \label{theorem:hardedge}
\end{theorem}

It is particularly interesting that the exponential attraction of
the zeros to the nodes of orthogonalization in $X_N$, that we have
seen is a feature of the asymptotics in subintervals of $[a,b]$
where the upper constraint is achieved by the equilibrium measure,
persists right up to the first and last nodes; in other words if
the upper constraint is achieved at $z=a$ then there is a zero
exponentially close to $x_{N,0}$ and if the upper constraint is
achieved at $z=b$ then there is a zero exponentially close to
$x_{N,N-1}$.

More is true, however.  From Proposition~\ref{prop:confine}, we
know that a zero $z_0$ of $\pi_{N,k}(z)$ that is exponentially
close to the first node $x_{N,0}$ must in fact satisfy the strict
inequality $z_0>x_{N,0}$.  Similarly, if a zero $z_0$ is
exponentially close to the last node $x_{N,N-1}$, then it must
satisfy the strict inequality $z_0<x_{N,N-1}$.  Going back to the
discussion in \S~\ref{sec:saturated}, we see that if there is a
hard edge at an endpoint of $[a,b]$, then in the saturated region
adjacent to the hard edge there can be no dislocations in the
pattern of Hurwitz zeros.  This is consistent with the fact (see
Theorem~\ref{theorem:upper}) that the function $\phi_\Gamma(z)$
does not vanish in any saturated region adjacent to endpoints of
$[a,b]$ so that there is no spurious zero.

\begin{remark}
The fact that the asymptotic formulae presented in
Theorem~\ref{theorem:hardedge} are in terms of the Euler gamma
function is directly related to the discrete nature of the
weights. In a sense, the poles of the functions
$\Gamma(1/2\pm\zeta_a)$ and $\Gamma(1/2\pm\zeta_b)$ are
``shadows'' of the poles of the Riemann-Hilbert problem that we
will discuss below in \S~\ref{sec;proof}.
\end{remark}

\section{Applications}\label{sec;applications}

\subsection{Discrete orthogonal polynomial
ensembles.}\label{sec;DOPensemble}

Recall that $X_N=\{ x_{N,n} \}_{n=0}^{N-1}$ is the set of nodes in
$(a,b)$.  In this section, we use the notation $w_N(x)$ for a
weight on $X_N$; to connect with our previous notation note simply
that for a node $x=x_{N,j}\in X_N$,
\begin{equation}
  w_N(x)= w_{N,j}\,.
\label{eq:weightrewrite}
\end{equation}
Consider the joint probability distribution of finding $k$
particles, say $P_1,\dots,P_k$, at respective positions
$x_1,\dots,x_k$ in $X_N$, to be given by the following expression:
\begin{equation}
\begin{array}{rcl}
\displaystyle \mathbb{P}(\mbox{particle $P_j$ lies at the site
$x_j$, for $j=1,\dots,k$})&=&
p^{(N,k)}(x_1,\dots,x_k)\\\\
&:=&\displaystyle \frac{1}{Z_{N,k}}\prod_{1\le i<j\le
k}(x_i-x_j)^2\prod_{j=1}^k w_N(x_j)\,,
\end{array}
\label{eq:prob}
\end{equation}
(we are using the symbol $\mathbb{P}(\mbox{event})$ to denote the
probability of an event) where $Z_{N,k}$ is a normalization
constant (or partition function) chosen so that
\begin{equation}
\sum_{\mbox{admissable configurations of $P_1,\dots,P_k$}}
p^{(N,k)}(x_1,\dots,x_k)=1\,. \label{eq:pnorm}
\end{equation}
Since the distribution function is symmetric in all $x_j$, we can
consider the particles $P_j$ to be either distinguishable or
indistinguishable, and only the normalization constant will depend
on this choice (the meaning of ``admissable configurations'' in
\eqref{eq:pnorm} is different in the two cases).
The statistical ensemble associated with the density function
\eqref{eq:prob} is called a {\em discrete orthogonal polynomial
ensemble}.

Discrete orthogonal polynomial ensembles arise in a number of
specific contexts (see for example, \cite{kurtj:shape,
Johansson01, Johansson00, BorodinO02}), with particular choices of
the weight function $w_N(\cdot)$ related (in cases we are aware
of) to classical discrete orthogonal polynomials.  It is of some
theoretical interest to determine properties of the ensembles that
are more or less independent of the particular choice of weight
function, at least within some class.  Such properties are said to
support the conjecture of {\em universality} within the class of
weight functions under consideration.

Some common properties of discrete orthogonal polynomial ensembles
can be read off immediately from the formula \eqref{eq:prob}.  For
example, the presence of the Vandermonde factor means that the
probability of finding two particles at the same site in $X_N$ is
zero.  Thus a discrete orthogonal polynomial ensemble always
describes an exclusion process.  This phenomenon is the discrete
analogue of the familiar {\em level repulsion} phenomenon in
random matrix theory. Also, since the weights are associated with
nodes, the interpretation is that configurations where particles
are concentrated in sets of nodes where the weight is larger are
more likely.

The goal of this section is to establish asymptotic formulae for
various statistics associated with the ensemble \eqref{eq:prob}
for a general class of weights in the continuum limit
$N\rightarrow\infty$ with the number of particles $k$ chosen so
that for some fixed rational $c\in (0,1)$, we have $k=cN$. Note
that the number of particles $k$ will have the same role as the
degree of orthogonal polynomials $k$. We use the same assumptions
on the nodes and weights as in the rest of the paper (see
\S~\ref{sec:C1}, \S~\ref{sec:C2}, and \S~\ref{sec:C3}). The main
idea is that, as is well-known, the formulae for all relevant
statistics of ensembles of the form \eqref{eq:prob} can be written
explicitly in terms of the discrete orthogonal polynomials
associated with the nodes $X_N$ and the weights
$w_{N,j}=w_N(x_{N,j})$.

To relate the statistics of interest to the discrete orthogonal
polynomials, we first define the so-called reproducing kernel
(Christoffel-Darboux kernel)
\begin{equation}
K_{N,k}(x,y):=\sqrt{w_N(x)w_N(y)}\sum_{n=0}^{k-1}p_{N,n}(x)p_{N,n}(y)\,,
\end{equation}
for $x,y$ in the nodes. Using the Christoffel-Darboux formula
\cite{Szego}, which holds for all orthogonal polynomials, even in
the discrete case, the sum on the right telescopes:
\begin{equation}
\begin{array}{rcl}
\displaystyle K_{N,k}(x,y)&=&\displaystyle \sqrt{w_N(x)w_N(y)}
\frac{c_{N,k-1}^{(k-1)}}{c_{N,k}^{(k)}}\cdot
\frac{p_{N,k}(x)p_{N,k-1}(y)-p_{N,k-1}(x)p_{N,k}(y)}{x-y} \\\\
&=&\displaystyle \sqrt{w_N(x)w_N(y)} \frac{\pi_{N,k}(x)\cdot
c_{N,k-1}^{(k-1)}p_{N,k-1}(y)-
c_{N,k-1}^{(k-1)}p_{N,k-1}(x)\cdot\pi_{N,k}(y)}{x-y}.
\end{array}
\end{equation}
Standard calculations (see, for example, \cite{Mehta91} or
\cite{TracyW93}) of random matrix theory, in the case of so-called
$\beta=2$ ensembles, yield the following exact formulae.  The
\emph{$m$-point correlation function} defined for $m\le k$ by
\begin{equation}
  R_m^{(N,k)}(x_1,\dots, x_m)
:= \frac{k!}{(k-m)!} \sum_{(x_{m+1}, \dots, x_k)\in X_N^{k-m}}
p^{(N,k)}(x_1, \dots, x_k)\,,
\end{equation}
where $X_N^p$ denotes the $p$th Cartesian power of $X_N$, can be
expressed in terms of the discrete orthogonal polynomials by the
formula
\begin{equation}
  R_m^{(N,k)}(x_1,\dots, x_m)
= \det \bigl( K_{N,k} (x_i, x_j) \bigr)_{1\le i,j\le m}\,.
\end{equation}
For any set $B\subset X_N$, the $1$-point correlation function has
the following interpretation:
\begin{equation}
  \sum_{x\in B} R^{(N,k)}_1(x)
= \mathbb{E}( \text{number of particles in $B$})\,,
\end{equation}
where $\mathbb{E}$ denotes the expected value.  Similarly, the
$2$-point correlation function has the following interpretation:
\begin{equation}
  \sum_{x,y\in B} R^{(N,k)}_2(x,y)
= \mathbb{E}( \text{number of (ordered) pairs of particles in
$B$}).
\end{equation}
Furthermore, the statistic defined for a set $B\subset X_N$ and
$m\le{\rm min}(\#B,k)$:
\begin{equation}\label{eq;Am1}
A_m^{(N,k)}(B) := \mathbb{P} ( \text{there are precisely $m$
particles in the set $B$})
\end{equation}
(this probability is automatically zero if $m>\#B$ by exclusion)
is well-known to be expressible by the exact formula
\begin{equation}\label{eq;Am2}
  A_m^{(N,k)}(B) =
\frac1{m!} \biggl(-\frac{d^m}{dt^m}\biggr)\biggl|_{t=1} \det
\bigl( 1-tK_{N,k}\bigl|_B \bigr)\,,
\end{equation}
where $K_{N,k}$ is the operator (in this case a finite matrix,
since $B$ is contained in the finite set $X_N$) acting in
$\ell^2(X_N)$ given by the kernel $K_{N,k}(x,y)$, and
$K_{N,k}\bigl|_B$ denotes the restriction of $K_{N,k}$ to
$\ell^2(B)$.

This is by no means an exhaustive list of statistics that can be
directly expressed in terms of the orthogonal polynomials
associated with the (discrete) weight $w_N(\cdot)$.  For example,
one may consider the fluctuations and in particular the variance
of the number of particles in an interval $B\subset X_N$.  The
continuum limit asymptotics for this statistic were computed in
\cite{Johansson00} for the Krawtchouk ensemble (see Proposition
2.5 of that paper) with the result that the fluctuations are
Gaussian; it would be of some interest to determine whether this
is special property of the Krawtchouk ensemble, or a universal
property of a large class of ensembles.  Also, there are
convenient formulae for statistics associated with the spacings
between particles; the reader can find such formulae in section
5.6 of the book \cite{deift}.


Depending on the location of interest, we have different results.
We distinguish again three regions: bands, voids and saturated
regions.

\subsubsection{In a band.}

\begin{theorem}[Universality of the discrete sine kernel]\label{theorem;bulk}
For a node $x\in X_N$ lying in a band $I$,
\begin{equation}
K_{N,k}(x,x) = \frac{c}{\rho^0(x)} \frac{d\mu^c_{\rm min}}{dx}(x)
\left( 1+
  O\left( \frac{\log(N)}{N} \right) \right)\,,
\end{equation}
where the error is uniform in compact subsets.  For distinct nodes
$x$ and $y$ in $I$,
\begin{equation}
  K_{N,k}(x,y) = \frac{O(1)}{N\cdot(x-y)},
\end{equation}
where $O(1)$ is uniform for $x$ and $y$ in a compact subsets. Also
with a given node $x\in I$, and for $\xi$ and $\eta$ such that
\begin{equation}
x+ \frac{\xi}{N\rho^0(x)K_{N,k}(x,x)}\in X_N\hspace{0.2
in}\mbox{and} \hspace{0.2 in} x+
\frac{\eta}{N\rho^0(x)K_{N,k}(x,x)}\in X_N\,,
\end{equation}
we have
\begin{equation}
\frac{1}{K_{N,k}(x,x)}
   K_{N,k}\left(x+ \frac{\xi}{N\rho^0(x)K_{N,k}(x,x)},
x+ \frac{\eta}{N\rho^0(x)K_{N,k}(x,x)}\right)= \frac{\sin(
\pi\cdot(\xi-\eta)
  )}{\pi\cdot(\xi-\eta)} +
  O\left( \frac{\log(N)}{N} \right)\,,
\end{equation}
where the error is uniform for $x$ in a compact subset of the band
$I$, and $\xi$ and $\eta$ in a compact set of $\mathbb{R}$.
\end{theorem}

\begin{remark}
Let $\psi^c_{\rm min}=d \mu^c_{\rm min}/dx$. By the same analysis,
we have the same limit for
\begin{equation}
  \frac{\rho^0(x)}{c\psi^c_{\rm min}(x)}
   K_{N,k}\left(x+ \frac{\xi}{Nc\psi^c_{\rm min}(x)},
x+ \frac{\eta}{Nc\psi^c_{\rm min}(x)}\right)
\end{equation}
with the same error bound.
\end{remark}

\begin{remark}
We believe that the logarithmic term $\log(N)$ in the error can be replaced by $1$ whenever
the nodes are
equally spaced ({\em cf.} the remark immediately following Theorem
\ref{theorem:band}), and it may be the case that such an improved estimate holds more
generally. In any case, the
improved factor of $1/N$ compared to either $1/N^{1/3}$ or $1/N^{2/3}$ as one might expect from the form of the error term $\varepsilon_N(z)$ in Theorems
\ref{theorem:band} and \ref{theorem:upper} is due to the particular structure
of the kernel $K_{N,k}$. Operators with this special type of kernel are
called {\em integrable operators} (see \cite{IIKS} and
\cite{deiftint}).
\end{remark}

Let the operator $\mathcal{S}_x$ act on $\ell^2(\mathbb{Z})$ with
the kernel
\begin{equation}
  \mathcal{S}_x(i,j)
  = \frac{\displaystyle\sin \left(\frac{\pi c\psi^c_{\rm min}(x)}{\rho^0(x)}
\cdot (i-j)\right)}{\pi \cdot (i-j)}\,,
  \qquad i,j\in\mathbb{Z}\,.
\end{equation}
Recall the formula \eqref{eq;Am2} for $A_m^{(N,k)}(B)$ and its
interpretation \eqref{eq;Am1} as a probability.
\begin{theorem}[Asymptotics of local occupation probabilities]
\label{theorem;bulk2} Let $B_N\subset X_N$ be a set of $M$ nodes
of the form
\begin{equation}
B_N=\{x_{N,j},x_{N,j+k_1},x_{N,j+k_2},\dots,x_{N,j+k_{M-1}}\}
\end{equation}
where $\#B_N=M$ is independent of $N$, and where
\begin{equation}
0<k_1<k_2<\cdots <k_{M-1} \hspace{0.2 in}\mbox{all in
$\mathbb{Z}$}
\end{equation}
are also all independent of $N$.  Set
$\mathbb{B}_N:=\{0,k_1,k_2,\dots,k_{M-1}\}\subset \mathbb{Z}$.
Suppose also that as $N\rightarrow\infty$, $x_{N,j}=\min
B_N\rightarrow x$ with $x$ lying in a band (and hence the same
holds for $x_{N,j+k_{M-1}}=\max B_N$).  Then, as
$N\rightarrow\infty$,
\begin{equation}\label{eq;limittoS1}
  \det \left(1-tK_{N,k}\bigl|_{B_N} \right)
  = \det\left( 1 - t\mathcal{S}_x\bigl|_{\mathbb{B}_N} \right) +
O\left(\frac{\log(N)}{N}\right)\,,
\end{equation}
for $t$ in a compact set in $\cx$, and
\begin{equation}\label{eq;limittoS2}
  A_m^{(N,k)}(B_N)
  = \frac1{m!}\left(-\frac{d^m}{dt^m}\right)\biggl|_{t=1}
  \det\left( 1 -t \mathcal{S}_x\bigl|_{\mathbb{B}_N} \right) +
  O\left(\frac{\log(N)}{N}\right)\,.
\end{equation}
\end{theorem}


\subsubsection{In voids and saturated regions.}

\begin{theorem}[Exponential asymptotics of the one-point function in voids]
Let $\Gamma$ be a void interval.  For each compact subset $F$ of
$\Gamma$, there is a constant $K_F>0$ such that
\begin{equation}
K_{N,k}(x,x) = O(e^{-K_FN})\hspace{0.2 in}\mbox{as}\hspace{0.2 in}
N\rightarrow\infty
\end{equation}
holds for all nodes $x\in X_N\cap F$. Also, for distinct nodes $x,
y$ in $X_N\cap F$ we have
\begin{equation}
 K_{N,k}(x,y) = \frac{O(e^{-K_FN})}{x-y}\,.
\end{equation}
\label{theorem:onepointvoid}
\end{theorem}
Thus, the one-point function is exponentially small in void
intervals as $N\rightarrow\infty$, going to zero with a decay rate
that is determined by the size of $\delta E_c/\delta\mu-\ell_c$ at
the node $x$.


\begin{theorem}[Exponential asymptotics of the one-point function in
saturated regions] Let $\Gamma$ be a saturated region.  For each
compact subset $F$ of $\Gamma$, there is a constant $K_F>0$ such
that
\begin{equation}
K_{N,k}(x,x)=1+O(e^{-K_FN})\hspace{0.2 in}\mbox{as}\hspace{0.2 in}
N\rightarrow\infty
\end{equation}
holds for all nodes $x\in X_N\cap F$. Also, for distinct nodes $x,
y$ in $X_N\cap F$ we have
\begin{equation}
 K_{N,k}(x,y) = \frac{O(e^{-K_FN})}{x-y}\,.
\end{equation}
\label{theorem:onepointsaturatedregion}
\end{theorem}
Therefore the one-point function is exponentially close to one in
saturated
regions.

\subsection{Random rhombus tiling of a hexagon.}\label{sec;tiling}

\begin{figure}[ht]
 \centerline{\epsfig{file=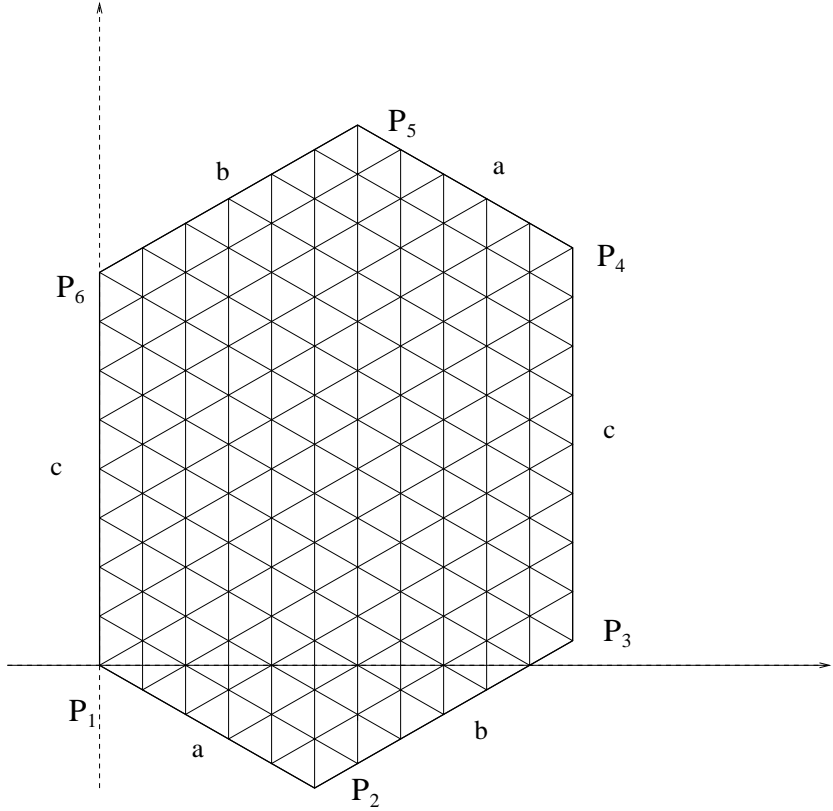, width=11cm}}
 \caption{\em
The $abc$-hexagon with vertices $P_1,\cdots, P_6$, and the lattice
$\mathcal{L}$} \label{fig-tilingbase}
\end{figure}

Let $a,b,c$ be positive integers, and consider a hexagon with
sides of lengths that proceed in counter-clockwise order,
$b,a,c,b,a,c$. All interior angles of this hexagon are equal and
measure $2\pi/3$ radians. We call this an {\em $abc$-hexagon}.
See Figure \ref{fig-tilingbase} for an example of an
$abc$-hexagon. We denote by $\mathcal{L}$ the lattice points
indicated in Figure \ref{fig-tilingbase}. By definition,
$\mathcal{L}$ includes the points on the sides $(P_6,P_1)$,
$(P_1,P_2)$, $(P_2,P_3)$, and $(P_3,P_4)$, but excludes the points
on the sides $(P_4,P_5)$ and $(P_5,P_6)$.

\begin{figure}[ht]
 \centerline{\epsfig{file=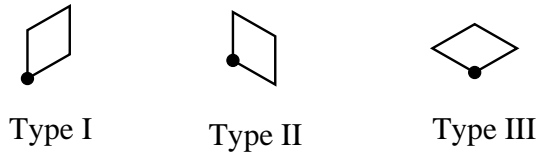, width=7cm}}
 \caption{\em The three types of rhombi; the position of each rhombus is
indicated with a dot.} \label{fig-tilingrhombi}
\end{figure}

\begin{figure}[ht]
 \centerline{\epsfig{file=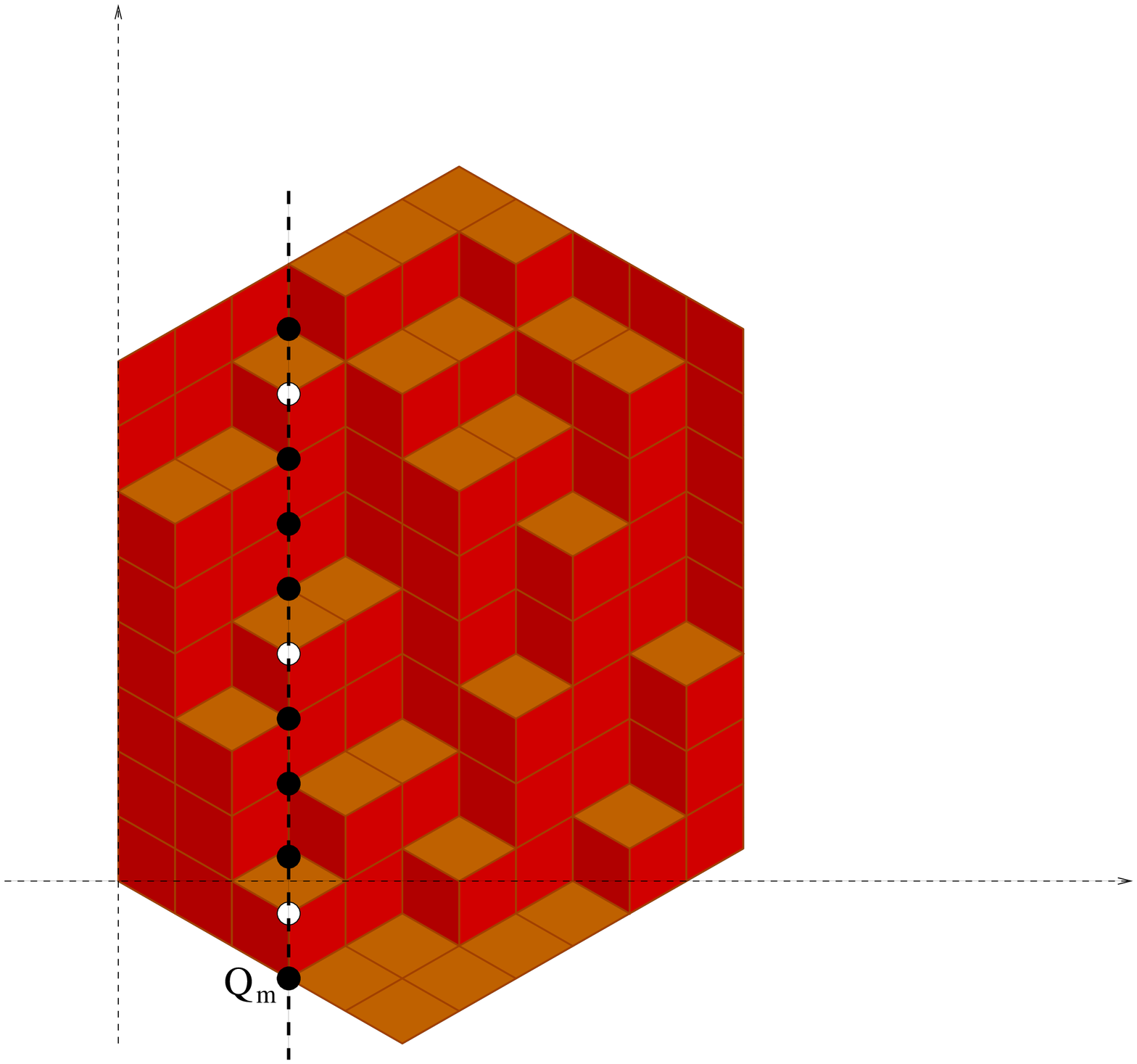, width=11cm}}
 \caption{\em A rhombus tiling of the $abc$-hexagon, and
the lattice $\mathcal{L}_m$ when $m=3$; holes are represented by
white dots and particles are represented by black dots.}
\label{fig-tilingandLm}
\end{figure}

Consider tiling the $abc$-hexagon with rhombi having sides of unit
length. Such rhombi come in three different types (orientations)
that we refer to as type I, type II, and type III; see
Figure~\ref{fig-tilingrhombi}.
Rhombi of types I and II are sometimes collectively called
\emph{horizontal rhombi}, while rhombi of type III are sometimes
called \emph{vertical rhombi}. The ``position" of each rhombus tile in
the hexagon is a specific lattice point in $\mathcal{L}$ defined as
indicated in Figure \ref{fig-tilingrhombi}. See
Figure~\ref{fig-tilingandLm} for an example of a rhombus tiling.

MacMahon's formula \cite{MacMahon60} gives the total number of all
possible rhombus tilings of the $abc$-hexagon as the expression
\begin{equation}
  \prod_{i=1}^a \prod_{j=1}^b \prod_{k=1}^c \frac{i+j+k-1}{i+j+k-2}\,.
\end{equation}
Consider the set of all rhombus tilings equipped with uniform
probability.  It is of some current interest to determine the
behavior of various corresponding statistics of this ensemble in
the limit as $a,b,c\to\infty$.

In the scaling limit of $n\to\infty$ where
\begin{equation}\label{eq;largehexagon}
  a= \alpha n\,, \qquad b= \beta n\,, \qquad c=\gamma n\,,
\end{equation}
with fixed $\alpha, \beta, \gamma >0$, the regions near the six
corners are ``frozen" or ``polar zones'', while the inside of the hexagon is ``temperate".
Cohn, Larsen
and Propp \cite{CohnLP98} showed that in such a limit, the
expected shape of the boundary of the frozen regions is given by
the inscribed ellipse. Moreover, the same authors also computed
the expected number of vertical rhombi in an arbitrary set $U\in
\R^2$. However, this calculation was provided without specific
error bounds.  In \cite{Johansson00}, Johansson proved a large
deviation result for the boundary shape, and also proved weak
convergence of the marginal probability of finding, say, a
vertical tile near a given location in a temperate region.   The same paper also contains
an investigation of the related Aztec diamond tiling model and a proof that the fluctuation of the boundary in this model is governed (in a proper scaling
limit) by the so-called Tracy-Widom law for large random matrices from
the Gaussian unitary ensemble
\cite{TW2}. The same is expected to be true for rhombus tilings of hexagons,
but this is still open.

In \cite{Johansson00}, Johansson expresses the
induced probability for a given configuration of vertical or horizontal rhombi on a given sublattice
in terms of discrete orthogonal polynomial ensembles with Hahn
or associated Hahn weights. Even though the Hahn weight is a
classical weight, the relevant asymptotics for Hahn polynomials have not been
previously established.  However, the asymptotics of the previous sections may now be applied to the
special case of the Hahn polynomials, and this yields new results
for the asymptotic properties of the hexagon tiling problem (see Theorems~\ref{theorem:hexagonstrong} and \ref{theorem:hexagonsine} below).


We first state the relation between hexagon tiling and discrete
orthogonal polynomial ensembles. We will assume without loss of
generality that $a\ge b$ (by the symmetry of the hexagon, the case
when $a\le b$ is completely analogous). Consider the $m^{th}$
vertical line of the lattice $\mathcal{L}$ counted from the left.
We denote by $\mathcal{L}_m$ the intersection of this line and
$\mathcal{L}$.  In a given tiling, the points in $\mathcal{L}_m$
correspond to positions (in the sense defined above) of a number
of rhombi of types I, II, and III.  We call the positions of
horizontal rhombi the \emph{particles}, and the positions of
vertical rhombi the \emph{holes}.  See
Figure~\ref{fig-tilingandLm} for an example of $\mathcal{L}_m$
when $m=3$, illustrating the corresponding particles and holes.

The uniform probability distribution on the ensemble of tilings
induces the probability distribution for finding particles and
holes at particular locations in the one-dimensional finite
lattice $\mathcal{L}_m$. A surprising result is due to Johansson
\cite{Johansson00} which states that the induced probability
distribution functions for holes and particles are both discrete
orthogonal polynomial ensembles with Hahn and associated Hahn
weight functions respectively (see \eqref{eq;Hahn1} and \eqref{eq;aHahn1}).

Let $Q_m$ be the lowest point in the sublattice $\mathcal{L}_m$.
On the sublattice $\mathcal{L}_m$, there are always $c$ particles,
and $L_m$ holes. We set $\gamma_m=c+L_m-1$. Now, let $x_1<\cdots <
x_c$, where $x_j\in\{ 0,1,2,\dots, \gamma_m\}$, denote the
(ordered) distances of the particles in $\mathcal{L}_m$ from
$Q_m$, and let $\xi_1< \cdots < \xi_{L_m}$, where $\xi_j \in\{
0,1,2,\dots, \gamma_m\}$, denote the distances of the holes in
$\mathcal{L}_m$ from $Q_m$. In particular, we then have $\{ x_1,
\dots, x_c\}\cup \{ \xi_1, \dots, \xi_{L_m}\} = \{0,1,2,\cdots,
\gamma_m\}$. Let $\tilde{P}_m(x_1,\dots, x_c)$ denote the
probability of finding the particle configuration $x_1, \cdots,
x_c$, and let $P_m(\xi_1,\dots, \xi_{L_m})$ denote the probability
of finding the hole configuration $\xi_1, \cdots, \xi_{L_m}$.

\begin{prop}[Theorem 4.1 of \cite{Johansson00}]
Let $a,b,c\ge 1$ be given integers with $a\ge b$. Set $a_m:=
|a-m|$ and $b_m:= |b-m|$. Then
\begin{equation}
  \tilde{P}_m(x_1,\dots, x_c)
 = \frac1{\tilde{Z}_m} \prod_{1\le j<k\le c} (x_j-x_k)^2
  \prod_{j=1}^{c} \tilde{w}(x_j)\,,
\end{equation}
where $\tilde{Z}_m$ is the normalization constant (partition
function), and where the weight function is the associated Hahn
weight
\begin{equation}
\tilde{w}(n):=
 w^{\rm Assoc}_{N,n}(a_m+1, b_m+1)
= \frac{\tilde{C}}{n!(a_m+n)!(N-n-1)!(N-n-1+b_m)!}\,,
\end{equation}
for a certain constant $\tilde{C}$. Also,
\begin{equation}
  P_m(\xi_1,\dots, \xi_{L_m}) =
  \frac1{Z_m} \prod_{1\le j<k\le L_m} (\xi_j-\xi_k)^2
  \prod_{j=1}^{L_m} w(\xi_j)\,,
\end{equation}
where $Z_m$ is the normalization constant, and where the weight
function is the Hahn weight
\begin{equation}
  w(n)
:= w^{\rm Hahn}_{N,n}(a_m+1, b_m+1) =
C\frac{(n+a_m)!(N-n-1+b_m)!}{n!(N-n-1)!}\,,
\end{equation}
for a certain constant $C$.
\end{prop}


With
\begin{equation}
  m=\tau n\,,
\end{equation}
for $\gamma>0$, the scaling \eqref{eq;largehexagon} is precisely
the same scaling that we analyzed in
\S~\ref{sec;DOPensemble}. Also we can explicitly compute the
equilibrium measure for Hahn and associated Hahn using either the
result of \cite{KuijlaarsR98}, or solving the variational problem
as in \cite{DKM}, which will appear in the full version of this
paper. The calculations of the equilibrium measure and the
one-point correlation function imply that as $n\to\infty$, the
one-dimensional lattice $\mathcal{L}_m$, after rescaling to finite
size independent of $n$, consists of three disjoint intervals: one
band, surrounded by two gaps (either saturated regions or voids,
depending on parameters). The saturated regions and voids
correspond to the frozen regions or polar zones, while the central band is a
temperate region. Hence in particular, the endpoints the
band when considered as functions of $\tau$ determine the typical
shape of the boundary between the polar and temperate zones of the
rescaled $abc$-hexagon. Also the one-point function converges
pointwise except at the band edges, or at the boundary, to the
equilibrium measure. This was conjectured in \cite{Johansson00}
(including the edges), in which weak convergence was obtained.
Moreover, our computation of the one-point correlation function
provides the relevant error bounds, when we consider sets $U$
contained in a single line $\mathcal{L}_m$. One expects that with
additional analysis of the same formulae it should be possible to
show that the error is locally uniform with respect to $\tau$, in
which case the same bounds should hold for more general regions
$U\in \R^2$. We state our result in this direction as follows.
\begin{theorem}[Strong asymptotics
with explicit error bounds] On the line $\mathcal{L}_m$, where
$m=\tau n$ and $\tau$ is fixed as $n\to\infty$, the scaled holes
$\xi_j/n$ lying in the polar zones, uniformly bounded away from
the rescaled expected boundary between the polar and temperate
zones, have a one-point function asymptotically convergent to
either $1$ (in the polar zones near the vertices $P_2$ and $P_5$),
or to $0$ (in the polar zones near the vertices $P_1$, $P_3$,
$P_4$, and $P_6$), with an exponential rate of convergence of the
order $O(e^{-Kn})$ for some constant $K>0$.  The one-point
correlation function for the scaled holes $\xi/n$ in the temperate
zone converges to the corresponding equilibrium measure with an
error of the order $O(1/n)$, which is uniformly valid away
from the rescaled expected boundary between the polar and
temperate zones. \label{theorem:hexagonstrong}
\end{theorem}

\begin{remark}
We give the above error estimate as $O(1/n)$ rather than $O(\log(n)/n)$ because
the Hahn and associated Hahn polynomials are orthogonal on a set of nodes $X_n$
that are equally spaced.
\end{remark}

In the temperate zone, in addition to the one-point function,
which is the marginal distribution, we can control all
$k$-point correlation functions under proper scaling. One such consequence is
the following theorem on the scaling limit for the locations of
the holes.
\begin{theorem}[Discrete sine kernel correlations]
Let $x>0$ be rational such that $nx\in \mathbb{Z}_N$ and such that
$nx$ is in the temperate zone away from the expected boundary
between the polar and temperate zones with uniform order in $n$.
Let $B_m= \{nx, nx+j_1, nx+j_2, \cdots, nx+j_M\}$, and set
$\mathbb{B}= \{0, j_1, j_2, \cdots, j_M\}$.  Then
\begin{equation}
  \lim_{n\to\infty}
\mathbb{P}(\text{there are precisely $p$ holes in the set $B_m$})
= \frac1{p!} \biggl( -\frac{d^p}{dt^p}\biggr)\biggl|_{t=1} \det
\left( 1-t \mathcal{S}|_{\mathbb{B}} \right),
\end{equation}
where $S$ acts on $\ell^2(\mathbb{Z})$ with the kernel
\begin{equation}
  S(i,j) = \frac{\sin ( c(x)(i-j))}{\pi (i-j)}\,,
\end{equation}
for some constant $c(x)$. \label{theorem:hexagonsine}
\end{theorem}

\begin{remark}
All of the results we have written down for holes have analogous
statements in terms of particles using the duality relation
between the Hahn and associated Hahn weights that will be
explained in \S~\ref{sec:reverse}.
\end{remark}

\begin{remark}
Once one obtains the asymptotics near the band edge of the
equilibrium measure for discrete orthogonal polynomial
ensembles, fluctuation statistics of the boundary curve will be computable.
It is conjectured in \cite{Johansson00} that the limiting law at the
band edge is the Tracy-Widom distribution known from the Gaussian unitary ensemble of
random matrix theory.
\end{remark}

\section{Riemann-Hilbert Problems for Discrete Orthogonal
Polynomials}\label{sec;proof}

In this section, we discuss the main ideas of asymptotic analysis
of discrete orthogonal polynomials via a Riemann-Hilbert problem.

\subsection{The fundamental Riemann-Hilbert problem.}
We first introduce the Riemann-Hilbert problem characterization of
discrete orthogonal polynomials. For $k\in\Z$ consider the matrix
$\mat{P}(z;N,k)$ solving the following problem, which is a
discrete version of the analogous problem for continuous weights
first used in \cite{FIK}.
\begin{rhp}
Find a $2\times 2$ matrix
$\mat{P}(z;N,k)$ with the following
properties:
\begin{enumerate}
\item
{\bf Analyticity}: $\mat{P}(z;N,k)$ is an analytic function of $z$ for
$z\in\cx\setminus X_N$.
\item
{\bf Normalization}: As $z\rightarrow\infty$,
\begin{equation}
\mat{P}(z;N,k)\left(\begin{array}{cc}z^{-k} & 0 \\ 0 & z^k\end{array}
\right)={\mathbb I} + O\left(\frac{1}{z}\right)\,.
\label{eq:norm}
\end{equation}
\item
{\bf Singularities}: At each node $x_{N,j}$, the first column of
$\mat{P}$ is analytic and the second column of $\mat{P}$ has a simple
pole, where the residue satisfies the condition
\begin{equation}
\mathop{\rm Res}_{z=x_{N,j}}\mat{P}(z;N,k)=\lim_{z\rightarrow x_{N,j}}
\mat{P}(z;N,k)\left(\begin{array}{cc}0 & w_{N,j}\\ 0 & 0\end{array}\right)
=\left(\begin{array}{cc}0 & w_{N,j}P_{11}(x_{N,j};N,k)\\
0 & w_{N,j}P_{21}(x_{N,j},N,k)\end{array}\right)
\label{eq:poles}
\end{equation}
for $j=0,\dots,N-1$.
\end{enumerate}
\label{rhp:DOP}
\end{rhp}
\begin{prop}
Riemann-Hilbert Problem~\ref{rhp:DOP} has a unique solution when $0\le
k\le N-1$.  In this case,
\begin{equation}
\mat{P}(z;N,k)=
\left(\begin{array}{cc}
\pi_{N,k}(z) & \displaystyle
\sum_{j=0}^{N-1}\frac{w_{N,j}\pi_{N,k}(x_{N,j})}{z-x_{N,j}}\\\\
c_{N,k-1}^{(k-1)}p_{N,k-1}(z) &\displaystyle \sum_{j=0}^{N-1}
\frac{w_{N,j}c_{N,k-1}^{(k-1)}p_{N,k-1}(x_{N,j})}{z-x_{N,j}}
\end{array}\right)
\label{eq:usoln}
\end{equation}
if $k>0$ and
\begin{equation}
\mat{P}(z;N,0)=
\left(\begin{array}{cc}
1 &\displaystyle \sum_{j=0}^{N-1}\frac{w_{N,j}}{z-x_{N,j}}\\\\
0 & 1\end{array}\right)\,.
\label{eq:usolnzero}
\end{equation}
\label{prop:solnrhp}
\end{prop}

We analyze this Riemann-Hilbert problem asymptotically as
$N,k\to\infty$ by adapting the Deift-Zhou procedure developed in
\cite{DZ1} and subsequent work. Due to the conditions on the pole
and the separation of the zeros of discrete orthogonal
polynomials, we have two difficulties which we mentioned in
\S~\ref{sec:method}. In the following two sections, we describe the
main techniques we developed to overcome these difficulties. A complete
asymptotic analysis of the above Riemann-Hilbert problem will
appear in the full version of this paper.

In order to apply the usual Deift-Zhou method, we will transform
the above Riemann-Hilbert problem into a Riemann-Hilbert problem
with jump conditions on continuous contours: a transformation from
$\mat{P} \mapsto \mat{R}$. For this new Riemann-Hilbert problem,
the formal limit of accumulation of nodes can be rigorously
justified. However, in addition to the continuum limit of
$N\to\infty$ ($N$ being the number of nodes), we simultaneously
take the large degree limit $k\to\infty$. In the analysis of
\cite{DKMVZstrong, DeiftKMVZ99}, a method for this limit is to
conjugate the Riemann-Hilbert problem with the so-called
$g$-function that is defined as a log transform of the equilibrium
measure. It was crucial in the analysis of \cite{DKMVZstrong,
DeiftKMVZ99} for the continuous orthogonal polynomials that the
equilibrium measure has only a lower constraint. Actually the
lower constraint yields an exponentially decaying factor. Hence
the upper constraint condition \eqref{eq:constraints} of the
equilibrium measure for discrete orthogonal polynomials generates
an exponentially growing factor. In order to replace an
exponentially growing term with an exponentially decaying term, we
introduce another transformation \emph{before} we map $\mat{P}$ to
$\mat{R}$: we will introduce an intermediate Riemann-Hilbert
problem for $\mat{Q}$ so that $\mat{P} \mapsto \mat{Q} \mapsto
\mat{R}$ as an exact sequence of transformations. In \S~\ref{sec:reverse}, we discuss the
transformation $\mat{P}\mapsto \mat{Q}$ of reversing the
triangularity of residue matrices, that will eventually work in
our favor turning exponentially growing terms into exponentially decaying terms. In \S~\ref{sec:continuum}, we discuss the transformation
$\mat{Q}\mapsto\mat{R}$ from a Riemann-Hilbert problem with residue
conditions to a Riemann-Hilbert problem with jumps on continuous
contours.


\subsection{Selectively reversing triangularity of residue
matrices.}\label{sec:reverse}

Riemann-Hilbert Problem~\ref{rhp:DOP} involves residue matrices
that are upper-triangular.  It will be advantageous in general to
modify the matrix $\mat{P}(z;N,k)$ in order to arrive at a new
Riemann-Hilbert problem in which we have selectively reversed the
triangularity of the residue matrices near certain individual
nodes $x_{N,j}$. Let $\del\subset {\mathbb Z}_N$ where ${\mathbb
Z}_N:=\{0,1,2,\dots,N-1\}$ and denote the number of elements in
$\del$ by $\#\del$ and the complementary set ${\mathbb
Z}_N\setminus\del$ by $\nab$.  We will reverse the triangularity
for those nodes $x_{N,j}$ for which $j\in\del$. Consider the
matrix $\mat{Q}(z;N,k)$ related to the solution $\mat{P}(z;N,k)$
of Riemann Hilbert Problem~\ref{rhp:DOP} as follows:
\begin{equation}
\mat{Q}(z;N,k):=\mat{P}(z;N,k)\left[
\prod_{n\in\del} (z-x_{N,n})\right]^{-\sigma_3}=
\mat{P}(z;N,k)\left(\begin{array}{cc} \displaystyle \prod_{n\in\del}
(z-x_{N,n})^{-1} & 0 \\\\ 0 &\displaystyle \prod_{n\in\del}(z-x_{N,n})
\end{array}\right)\,.
\label{eq:PtoQ}
\end{equation}
It is direct to check that the matrix $\mat{Q}(z;N,k)$ is, for
$k\in{\mathbb Z}_N$, the unique solution of the following
Riemann-Hilbert problem.
\begin{rhp}
Given a subset $\del$ of ${\mathbb Z}_N$ of cardinality
$\#\del$, find a $2\times 2$ matrix $\mat{Q}(z;N,k)$ with the
following properties:
\begin{enumerate}
\item
{\bf Analyticity}: $\mat{Q}(z;N,k)$ is an analytic function of $z$ for
$z\in\cx\setminus X_N$.
\item
{\bf Normalization}: As $z\rightarrow\infty$,
\begin{equation}
\mat{Q}(z;N,k)\left(\begin{array}{cc}z^{\#\del-k} & 0 \\ 0 &
z^{k-\#\del}\end{array} \right)={\mathbb I} + O\left(\frac{1}{z}\right)\,.
\label{eq:norm-S}
\end{equation}
\item
{\bf Singularities}: At each node $x_{N,j}$, the matrix $\mat{Q}$ has
a simple pole.  If $j\in\nab$ where
$\nab:={\mathbb Z}_N\setminus\del$, then the first column
is analytic at $x_{N,j}$ and the pole is in the second column such
that the residue satisfies the condition
\begin{equation}
\mathop{\rm Res}_{z=x_{N,j}}\mat{Q}(z;N,k)=\lim_{z\rightarrow x_{N,j}}
\mat{Q}(z;N,k)\left(\begin{array}{cc}0 & \displaystyle w_{N,j}
\prod_{n\in\del}(x_{N,j}-x_{N,n})^2\\\\ 0 & 0\end{array}\right)
\label{eq:polesinS}
\end{equation}
for $j\in\nab$.
If $j\in\del$, then the second column is analytic at $x_{N,j}$ and the
pole is in the first column such that the residue satisfies the condition
\begin{equation}
\mathop{\rm Res}_{z=x_{N,j}}\mat{Q}(z;N,k)=\lim_{z\rightarrow x_{N,j}}
\mat{Q}(z;N,k)\left(\begin{array}{cc}0 & 0 \\\\
\displaystyle
\frac{1}{w_{N,j}}\mathop{\prod_{n\in\del}}_{n\neq j}(x_{N,j}-x_{N,n})^{-2}
& 0\end{array}\right)
\label{eq:polesnotinS}
\end{equation}
for $j\in\del$.
\end{enumerate}
\label{rhp:DOP-S}
\end{rhp}

\begin{figure}[h]
\begin{center}
\input{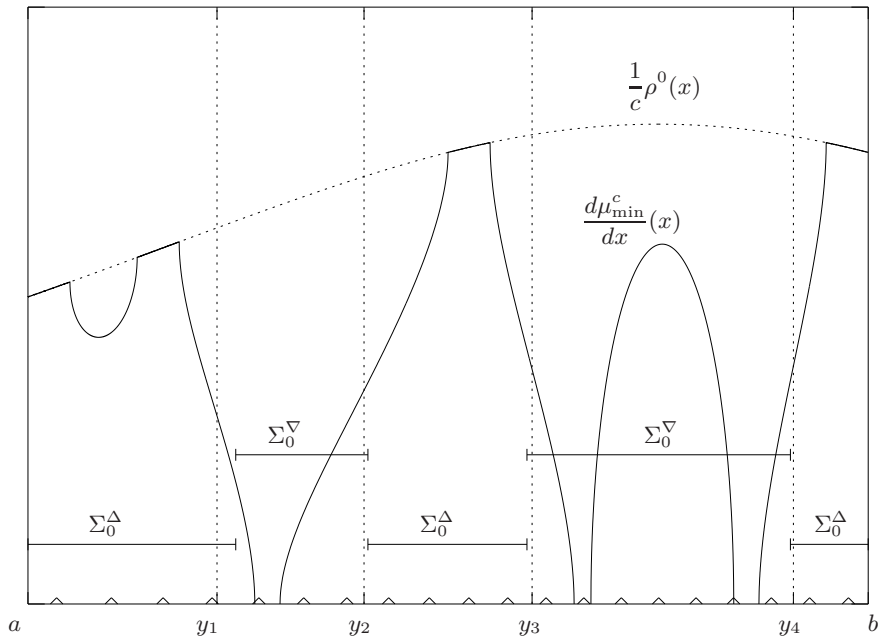}
\end{center}
\caption{\em A schematic diagram showing the relation of the
minimizer $\mu^c_{\rm min}(x)$ to the interval systems
$\Sigma_0^\nab$ and $\Sigma_0^\del$.  The nodes $x_{N,j}\in (a,b)$
are indicated on the $x$-axis with triangles; their density is
proportional to the upper constraint.  The common endpoints of subintervals of $\Sigma_0^\nab$ and
$\Sigma_0^\del$ converge as $N\rightarrow\infty$ to the points
$y_k$ indicated on the $x$-axis.
}
\label{fig:sigmaAB}
\end{figure}

Note that the $(21)$-entry of the residue matrix in
\eqref{eq:polesnotinS} is the reciprocal of the $(12)$-entry of the
residue matrix in \eqref{eq:polesinS}. When we make the choice
that $\del$ contains the saturated regions, the effect will be to turn
exponentially growing factors into exponentially decaying factors.
Let us be more specific about how we choose $\del$.  In each band $I_k$ lying between
a void and a saturated region we choose a point $y_k$, and ``quantize'' these to the lattice
$X_N$ by associating with each point a sequence $\{y_{k,N}\}_{N=0}^\infty$ converging
to $y_k$ as $N\rightarrow\infty$ with elements given by
\begin{equation}
N\int_a^{y_{k,N}}\rho^0(x)\,dx = \left\lceil N\int_a^{y_k}\rho^0(x)\,dx
\right\rceil
\end{equation}
where $\lceil u\rceil$ denotes the least integer greater than or equal to $u$.  Thus
$y_{k,N}$ lies asymptotically halfway between two consecutive nodes.  For each $N$,
these points are the common endpoints of two complementary systems of subintervals
of $(a,b)$.  We denote the union of open subintervals delineated by these points and containing no saturated regions by $\Sigma_0^\nab$.  The complementary system of
subintervals contains no voids and is denoted by $\Sigma_0^\del$.  See Figure~\ref{fig:sigmaAB}.


Based on this partitioning of $(a,b)$, the specific choice we make is that $\del$ is the set of indices $j\in\mathbb{Z}_N$ such that
$x_{N,j}\in\Sigma_0^\del$.

\begin{remark}
After the completion of this work, we learned that an analogous
transformation was used in \cite{BorodinO02} for a somewhat different
asymptotic analysis of a special choice of discrete orthogonal
polynomials (Askey-Lesky weights, or general Hahn weights).  In
\cite{BorodinO02}, the choice of the region $\Delta$ to reverse the
triangularity was made using intuition from representation theory.
(We thank A. Borodin for bringing this to our attention.)
For the limit of interest in this paper, we determine
the region $\del$ in order to reverse the triangularity in saturated
regions of the equilibrium measure while preserving triangularity in all voids. 
\end{remark}

\subsubsection{Dual families of discrete orthogonal polynomials.}

 The relation between Riemann-Hilbert
Problem~\ref{rhp:DOP} and Riemann-Hilbert Problem~\ref{rhp:DOP-S}
gives rise in a special case to a remarkable duality between pairs
of weights $\{w_{N,j}\}$ defined on the same set of nodes and
their corresponding families of discrete orthogonal polynomials
that comes up in applications.  Given nodes $X_N$ and weights
$\{w_{N,j}\}$, take $\del={\mathbb Z}_N$ and let
\begin{equation}
\overline{\mat{P}}(z;N,\overline{k}):=
\sigma_1 \mat{Q}(z;N,k) \sigma_1\,,\hspace{0.2 in}\mbox{where}
\hspace{0.2 in}\overline{k}:=N-k\,.
\end{equation}
Thus, we are reversing the triangularity at all of the nodes, and
swapping rows and columns of the resulting matrix.  It is easy to check
that $\overline{\mat{P}}(z;N,\overline{k})$ satisfies
\begin{equation}
\overline{\mat{P}}(z;N,\overline{k})\left(\begin{array}{cc}
z^{-\overline{k}} & 0 \\ 0 & z^{\overline{k}}\end{array}\right) =
{\mathbb I} + O\left(\frac{1}{z}\right)\hspace{0.2 in}\mbox{as}
\hspace{0.2 in}z\rightarrow\infty
\end{equation}
and is a matrix with simple poles in the second column at all nodes,
such that
\begin{equation}
\mathop{\rm Res}_{z=x_{N,j}}\overline{\mat{P}}(z;N,\overline{k})=
\lim_{z\rightarrow x_{N,j}}
\overline{\mat{P}}(z;N,\overline{k})
\left(\begin{array}{cc}0 & \overline{w}_{N,j}\\ 0 & 0\end{array}\right)
\end{equation}
holds for $j\in {\mathbb Z}_N$, where the ``dual weights''
$\{\overline{w}_{N,j}\}$ are defined by the identity
\begin{equation}
w_{N,j}\overline{w}_{N,j}\mathop{\prod_{n=0}^{N-1}}_{n\neq j}(x_{N,j}-x_{N,n})^2 = 1\,.
\label{eq:dualweightsdefine}
\end{equation}
Comparing with Riemann-Hilbert Problem~\ref{rhp:DOP} we see that
$\overline{P}_{11}(z;N,\overline{k})$ is the monic orthogonal
polynomial $\overline{\pi}_{N,\overline{k}}(z)$ of degree
$\overline{k}$ associated with the dual weights
$\{\overline{w}_{N,j}\}$.  In this sense, families of discrete
orthogonal polynomials always come in dual pairs.  An explicit
relation between the dual polynomials comes from the representation of
$\mat{P}(z;N,k)$ given by Proposition~\ref{prop:solnrhp}:
\begin{equation}
\begin{array}{rcl}
\displaystyle
\overline{\pi}_{N,\overline{k}}(z)&=&\displaystyle
\overline{P}_{11}(z;N,\overline{k})\\\\
&=&\displaystyle
P_{22}(z;N,k)\prod_{n=0}^{N-1}(z-x_{N,n})\\\\
&=&\displaystyle
\sum_{j=0}^{N-1}w_{N,j}\left[c_{N,k-1}^{(k-1)}\right]^2\pi_{N,k-1}(x_{N,j})
\mathop{\prod_{n=0}^{N-1}}_{n\neq j}(z-x_{N,n})\,.
\end{array}
\label{eq:explicitdual}
\end{equation}
Since the left-hand side is a monic polynomial of degree
$\overline{k}=N-k$ and the right-hand side is apparently a polynomial
of degree $N-1$, equation (\ref{eq:explicitdual}) furnishes $k$
relations among the weights and the normalization constants
$c_{N,k}^{(k)}$.

In particular, if we evaluate (\ref{eq:explicitdual}) for $z=x_{N,l}$
for some $l\in {\mathbb Z}_N$, then only one term from the sum on the
right-hand side survives and we find
\begin{equation}
\overline{\pi}_{N,\overline{k}}(x_{N,l})=\left[c_{N,k-1}^{(k-1)}\right]^2
w_{N,l}\mathop{\prod_{n=0}^{N-1}}_{n\neq l}
(x_{N,l}-x_{N,n})\cdot\pi_{N,k-1}(x_{N,l})\,,
\label{eq:Borodin}
\end{equation}
an identity relating values of each discrete orthogonal polynomial and
a corresponding dual polynomial at any given node.  The identity
\eqref{eq:Borodin} has also been derived by Borodin \cite{Borodin01}.

\begin{remark}
We want to point out that the notion of duality described here is
different from that explained in \cite{NikiforovSU91}.  The latter
generally involves relationships between families of discrete
orthogonal polynomials with two different sets of nodes of
orthogonalization.  For example, the Hahn polynomials are orthogonal
on a lattice of equally spaced points, and the polynomials dual to
the Hahn polynomials by the scheme of \cite{NikiforovSU91} are
orthogonal on a quadratic lattice for which $x_{N,n}-x_{N,n-1}$ is
proportional to $n$.  However, the polynomials dual to the Hahn
polynomials under the scheme described above are the associated Hahn
polynomials, which are orthogonal on the same equally-spaced nodes
as are the Hahn polynomials themselves.  The notion of duality we use
in this paper coincides with that described in
\cite{Borodin01} and is also equivalent to the ``hole/particle
transformation'' considered by Johansson \cite{Johansson00}.
\end{remark}

\subsection{Removal of poles in favor of jumps on
contours.}\label{sec:continuum}

\begin{figure}[h]
\begin{center}
\input{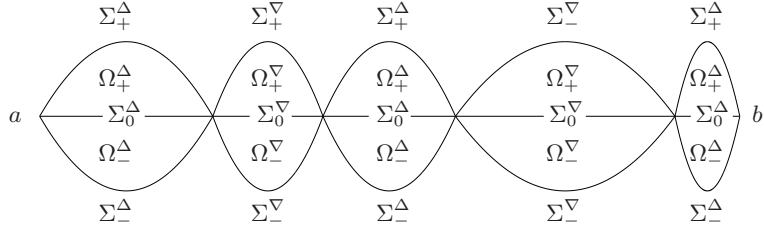}
\end{center}
\caption{\em The contour $\Sigma$ consists of the subintervals
$\Sigma_0^\nab$ and $\Sigma_0^\del$ as in
Figure~\ref{fig:sigmaAB} and associated contour segments
$\Sigma_+^\nab$ and $\Sigma_+^\del$ in the upper half-plane, and
$\Sigma_-^\nab$ and $\Sigma_-^\del$ in the lower half-plane. The
enclosed regions $\Omega_\pm^\nab$ and $\Omega_\pm^\del$ are also
indicated.  The contour $\Sigma$ lies entirely in the region of
analyticity of $V(x)$ and $\rho^0(x)$.  All components of $\Sigma$
are taken to be oriented from left to right.}
\label{fig:ABcontours}
\end{figure}

The deformations in this section
are based on similar ones first introduced by one of the authors
in \cite{Miller02}.  Let the analytic functions $\beta_\pm(z)$ be
given by
\begin{equation}
\beta_\pm(z):=\pm i\exp\left(\mp i\pi
N\int_z^b\rho^0(s)\,ds\right)\,. \label{eq:betadef}
\end{equation}
Note that by definition,
$\beta_+(x_{N,j})=\beta_-(x_{N,j})=(-1)^{N-1-j}$ for all
$N\in\nat$ and $j\in{\mathbb Z}_N$. Consider the contour $\Sigma$
illustrated in Figure~\ref{fig:ABcontours}.
>From the solution of Riemann-Hilbert Problem \ref{rhp:DOP-S} we
define a new matrix $\mat{R}(z)$ as follows.  Set
\begin{equation}
\mat{R}(z):=\mat{Q}(z;N,k)\left(\begin{array}{cc} 1 &
\displaystyle -\beta_\pm(z)e^{-NV_N(z)}\frac{\displaystyle
\prod_{j\in\del}(z-x_{N,j})}{\displaystyle\prod_{j\in\nab}
(z-x_{N,j})}\\\\
0 & 1\end{array}\right)\hspace{0.2 in}\mbox{for}\hspace{0.2 in}
z\in \Omega_\pm^\nab\,, \label{eq:QtoR1}
\end{equation}
\begin{equation}
\mat{R}(z):=\mat{Q}(z;N,k)\left(\begin{array}{cc}
1 & 0\\\\
\displaystyle -\beta_\pm(z)e^{NV_N(z)}\frac{\displaystyle
\prod_{j\in\nab}(z-x_{N,j})}
{\displaystyle\prod_{j\in\del}(z-x_{N,j})} & 1
\end{array}\right)\hspace{0.2 in}\mbox{for}\hspace{0.2 in}
z\in \Omega_\pm^\del\,, \label{eq:QtoR2}
\end{equation}
and for all other $z$ set $\mat{R}(z):=\mat{Q}(z;N,k)$.

The matrix $\mat{R}(z)$ is, for arbitrary $N\in\nat$ and
$k\in{\mathbb Z}_N$, the unique solution of the following
Riemann-Hilbert problem.
\begin{rhp}
Find a $2\times 2$ matrix $\mat{R}(z)$ with the following
properties:
\begin{enumerate}
\item
{\bf Analyticity}: $\mat{R}(z)$ is an analytic function of $z$ for
$z\in\cx\setminus \Sigma$.
\item
{\bf Normalization}: As $z\rightarrow\infty$,
\begin{equation}
\mat{R}(z)\left(\begin{array}{cc}z^{\#\del-k} & 0 \\ 0 &
z^{k-\#\del}\end{array} \right)={\mathbb I} +
O\left(\frac{1}{z}\right)\,.
\end{equation}
\item
{\bf Jump Conditions}: $\mat{R}(z)$ takes continuous boundary
values on $\Sigma$ from each connected component of
$\cx\setminus\Sigma$. Denoting the boundary values taken on the
left (right) by $\mat{R}_+(z)$ ($\mat{R}_-(z)$), we have
\begin{equation}
\begin{array}{rcl}
\displaystyle \mat{R}_+(z)&=&\displaystyle
\mat{R}_-(z)\left(\begin{array}{cc} 1 & \displaystyle
\pm\beta_\pm(z)e^{-NV_N(z)}
\frac{\displaystyle\prod_{j\in\del}(z-x_{N,j})}{\displaystyle
\prod_{j\in\nab}(z-x_{N,j})}\\\\
0 & 1\end{array}\right)\\\\
&&\displaystyle\hspace{0.2 in}\mbox{for}\hspace{0.2 in} z\in
\Sigma_\pm^\nab\,,
\end{array}
\end{equation}
\begin{equation}
\begin{array}{rcl}
\displaystyle\mat{R}_+(z)&=&\displaystyle
\mat{R}_-(z)\left(\begin{array}{cc}
1 & 0\\\\
\displaystyle \pm\beta_\pm(z)e^{NV_N(z)}
\frac{\displaystyle\prod_{j\in\nab}(z-x_{N,j})}
{\displaystyle\prod_{j\in\del}(z-x_{N,j})} &
1\end{array}\right)\\\\
&&\displaystyle\hspace{0.2 in}\mbox{for}\hspace{0.2 in} z\in
\Sigma_\pm^\del\,,
\end{array}
\end{equation}
\begin{equation}
\begin{array}{rcl}
\displaystyle\mat{R}_+(z)&=&\displaystyle
\mat{R}_-(z)\left(\begin{array}{cc} 1 & \displaystyle
(\beta_-(z)-\beta_+(z))e^{-NV_N(z)}
\frac{\displaystyle\prod_{j\in\del}(z-x_{N,j})}
{\displaystyle\prod_{j\in\nab}(z-x_{N,j})}\\\\
0 & 1\end{array}\right)\\\\
&&\displaystyle\hspace{0.2 in}\mbox{for}\hspace{0.2 in} z\in
\Sigma_0^\nab\,,\hspace{0.2 in}\mbox{and,}
\end{array}
\end{equation}
\begin{equation}
\begin{array}{rcl}
\displaystyle\mat{R}_+(z)&=&\displaystyle
\mat{R}_-(z)\left(\begin{array}{cc}
1 & 0\\\\
\displaystyle (\beta_-(z)-\beta_+(z))
e^{NV_N(z)}\frac{\displaystyle\prod_{j\in\nab}(z-x_{N,j})}
{\displaystyle\prod_{j\in\del}(z-x_{N,j})}
 & 1\end{array}\right)\\\\
&&\displaystyle\hspace{0.2 in}\mbox{for}\hspace{0.2 in} z\in
\Sigma_0^\del\,.
\end{array}
\end{equation}
Note that all off-diagonal entries of the jump matrices are
analytic nonvanishing functions on their respective contours.
\end{enumerate}
\label{rhp:R}
\end{rhp}
The significance of passing from Riemann-Hilbert
Problem~\ref{rhp:DOP-S} to Riemann-Hilbert Problem~\ref{rhp:R} is
that all poles have completely disappeared from the problem. All
boundary values of $\mat{R}(z)$ and the corresponding jump matrices
relating them are analytic functions. This means that
Riemann-Hilbert Problem~\ref{rhp:R} is sufficiently similar to
that introduced in \cite{FIK} for the continuous weight case that
it may, in principle, be analyzed by methods like those used in
\cite{DKMVZstrong,DeiftKMVZ99}.  The main obstruction at this
point is that the off-diagonal elements of the jump matrices for
$\mat{R}(z)$ are not exactly of the form $e^{NW(z)}$ for some
$W(z)$. This is a consequence of the fact that the sequence of transformations
$\mat{P}\mapsto\mat{Q}\mapsto\mat{R}$ is \emph{exact}, and one may observe
at this point that the desired form $e^{NW(z)}$ can be achieved by carefully taking a
natural continuum limit based on the assumptions on the nodes and weights set out
at the beginning of this announcement.  In other words, while the jump matrix for
$\mat{R}$ does not have the desired form, one may introduce an approximate Riemann-Hilbert problem for a matrix $\dot{\mat{R}}(z)$ for which the jump matrix indeed
has the desired form; part of the analysis then becomes the task of showing that $\mat{R}(z)$ and $\dot{\mat{R}}(z)$ are ``close''.

\begin{remark}
The difficulty that prevents us from obtaining asymptotic results near the band edges
can be traced back
precisely to the fact that the off-diagonal elements of the jump
matrices for $\mat{R}(z)$ are not exactly of the form $e^{NW(z)}$, but only
approximately so. In the local analysis, the (small) discrepancy between
the approximate and the exact form prevents us from obtaining the necessary
asymptotics up to the boundary of the contours. We hope to be able
to overcome this problem in a future work.
\end{remark}

The final aspect of our analysis that we would like to briefly describe is our choice of a $g$-function with which we stabilize the (approximate) Riemann-Hilbert problem for $\dot{\mat{R}}(z)$.  We introduce a new matrix unknown by the transformation
$\dot{\mat{S}}(z):=\dot{\mat{R}}(z)e^{(\#\del-k)g(z)\sigma_3}$ where
the $g$-function is given by
\begin{equation}
g(z)=\int_a^b\log(z-x)\rho(x)\,dx
\end{equation}
with density determined differently in the two types of intervals $\Sigma_0^\nab$ and
$\Sigma_0^\del$:
\begin{equation}
\rho(x):=\left\{
\begin{array}{ll}
\displaystyle\frac{c}{c-d}\frac{d\mu_{\rm min}^c}{dx}(x)\,, &
x\in\Sigma_0^\nab\\\\
\displaystyle\frac{c}{c-d}\left(\frac{d\mu_{\rm min}^c}{dx}(x)-\frac{1}{c}
\rho^0(x)\right)\,, & x\in\Sigma_0^\del\,.
\end{array}\right.
\end{equation}
Here of course $\mu_{\rm min}^c$ is the equilibrium measure, $c=k/N$, and
\begin{equation}
d:=\int_{\Sigma_0^\del}\rho^0(x)\,dx = \frac{\#\del}{N}\,.
\end{equation}
With this choice of $g(z)$, in conjunction with the choice of the set $\del$ described
above, the jump matrices for $\dot{\mat{S}}(z)$ are precisely of the type for which
the steepest-descent factorization technique can be applied.  Also, $\dot{\mat{S}}(z)$ is now normalized to the identity matrix for large $z$; the power asymptotics have been removed.  The complete details of the subsequent analysis, including rigorous error estimates,
will appear in the full version of the paper corresponding to this announcement.


\medskip
\noindent {\bf Acknowledgments.} We would like to thank Percy
Deift and Kurt Johansson for their interest and useful
conversations. Special thanks is due to Alexei Borodin who
informed us of his paper with Olshanski \cite{BorodinO02}. J. Baik
would like to thank the Institute for Advanced Study where a part
of work is done. K. McLaughlin wishes to thank T. Paul, F. Golse, and
the staff of the \'{E}cole Normal Superieur, Paris for their kind
hospitality. The research of J. Baik is supported in part by the National
Science Foundation under grant DMS-0208577. The research of K. McLaughlin
is supported in part by the National Science Foundation under
grants DMS-9970328 and DMS-0200749. The research of P. Miller is supported
in part by the National Science Foundation under grant
DMS-0103909. 

\bibliographystyle{plain}
\bibliography{ref}

\end{document}